\documentclass[a4paper]{article}

\usepackage{geometry}
\usepackage[utf8]{inputenc}
\usepackage[LGR,T1]{fontenc}
\usepackage[english]{babel}
\usepackage{textalpha}
\usepackage{lmodern}

\usepackage{amsmath}
\usepackage{amssymb}
\usepackage{amsthm}
\usepackage{bm}
\usepackage{graphicx}
\usepackage[svgnames]{xcolor}
\usepackage[babel]{csquotes}
\usepackage{titlecaps}
\Addlcwords{a an the and or nor but as at by for in of on per to vs}
\usepackage{doi}
\newcommand{\email}[1]{\href{mailto:#1}{\texttt{#1}}}
\usepackage{hyperref}
\addto\extrasenglish{%
}

\usepackage{thmtools}

\DeclareMathOperator{\e}{e}
\newcommand{\R}{\ensuremath{\mathbb{R}}}
\newcommand{\C}{\ensuremath{\mathbb{C}}}
\newcommand{\Z}{\ensuremath{\mathbb{Z}}}
\newcommand{\Y}{\ensuremath{\mathbb{H}}}

\newcommand{\W}{\ensuremath{\mathcal{W}}}
\newcommand{\N}{\ensuremath{\mathcal{N}}}
\renewcommand{\d}{\mathrm{d}}

\newcommand{\dy}{\d y}

\newcommand{\del}{\partial}

\newcommand{\ra}{\rightarrow}
\newcommand\ep{\varnothing}

\newcommand\p{\partial}
\newcommand\mm{\text}

\newcommand\si{\ensuremath{\mathcal{S}}}

\DeclareMathOperator{\diam}{diam}
\DeclareMathOperator{\vol}{Vol}

\DeclareMathOperator{\scal}{Scal}

\DeclareMathOperator{\length}{length}

\renewcommand{\epsilon}{\varepsilon}
\renewcommand{\Gamma}{\varGamma}
\renewcommand{\Phi}{\varPhi}
\renewcommand{\Omega}{\varOmega} % italic Omega

\newcommand{\diff}[3][]{\frac{\d^{#1}#2}{\d {#3}^{#1}}}

\renewcommand{\le}{\leqslant}
\renewcommand{\ge}{\geqslant}
\newcommand{\delG}{\del_{\mathrm{G}}} % Gromov Boundary
\newcommand{\delM}{\del_{\mathrm{M}}} % Martin Boundary
\newcommand{\F}{\mathcal{F}}

\DeclareMathOperator{\dist}{dist}
\newcommand{\eucl}{\ensuremath{\mathrm{Eucl}}}
\newcommand{\dd}{\textit{\dj\hspace{1pt}}}
\newcommand{\tdd}{\textit{\dh\hspace{1pt}}}
\newcommand{\red}{\mathcal{R}}
\newcommand{\E}{\mathcal{E}}
\newcommand\loc{\mathrm{loc}}
\renewcommand\c{\mathrm{c}}
\newcommand\Floc{\F_\loc}
\newcommand\Fc{\F_\c}
\renewcommand\L{\mathcal{L}}

\newcommand{\ap}{{\langle A \rangle}}
\renewcommand{\S}{\ensuremath{\mathcal{S}}}

\newcommand{\mule}{\preccurlyeq}

\newcommand{\thmref}[1]{\hyperref[#1]{\autoref*{#1} \enquote{\nameref*{#1}}}}
\newcommand{\thmnameref}[1]{\hyperref[#1]{\nameref*{#1} (\autoref*{#1})}}
\newcommand{\thmnonameref}[2]{\hyperref[#1]{#2 \ref*{#1}}}

\newtheoremstyle{note}% name
    {3pt}       % Space above
    {3pt}       % Space below
    {}          % Body font
    {}          % Indent amount
    {\itshape}  % Theorem head font
    {:}         % Punctuation after theorem head
    {.5em}      % Space after theorem head
    {}          % Theorem head spec (can be left empty, meaning ‘normal’)

\newtheoremstyle{plaincaps}% name is set in title case
    {}{}%
    {\itshape}{}%
    {\bfseries}{}%
    { }%
    {\thmname{#1}\thmnumber{ #2}\thmnote{ (\titlecap{#3})}}
\newtheoremstyle{plainnocaps}% name is not set in title case
    {}{}%
    {\itshape}{}%
    {\bfseries}{}%
    { }%
    {\thmname{#1}\thmnumber{ #2}\thmnote{ (#3)}}
\newtheoremstyle{definitioncaps}% name is set in title case
    {}{}%
    {}{}%
    {\bfseries}{}%
    { }%
    {\thmname{#1}\thmnumber{ #2}\thmnote{ (\titlecap{#3})}}
\newtheoremstyle{definitionnocaps}% name is not set in title case
    {}{}%
    {}{}%
    {\bfseries}{}%
    { }%
    {\thmname{#1}\thmnumber{ #2}\thmnote{ (#3)}}

\theoremstyle{plaincaps}

\newtheorem{theorem}{Theorem}[section]
\newtheorem{proposition}[theorem]{Proposition}
\newtheorem{defprop}[theorem]{Definition/Proposition}
\newtheorem{lemma}[theorem]{Lemma}
\newtheorem{corollary}[theorem]{Corollary}

\theoremstyle{plainnocaps}

\theoremstyle{definitioncaps}

\newtheorem{definition}[theorem]{Definition}
\newtheorem{property}[theorem]{Property}

\theoremstyle{definitionnocaps}

\theoremstyle{remark}
\newtheorem{remark}[theorem]{Remark}
\newtheorem{example}[theorem]{Example}
\newtheorem{remarks}[theorem]{Remarks}
\newtheorem{examples}[theorem]{Examples}

\hypersetup{%
	bookmarksopen,
	bookmarksnumbered,
	linkcolor=DarkGreen,
	citecolor=Navy,
	urlcolor=Navy,
	pdftitle={Potential Theory on Gromov Hyperbolic Spaces},
	pdfauthor={Matthias Kemper, Joachim Lohkamp}
}

\begin{document}

\title{Potential Theory on Gromov Hyperbolic Spaces}
\author{Matthias Kemper \and Joachim Lohkamp}
%\date{\vspace*{-0.9cm}}

%\institute{M. Kemper \at
%              Mathematisches Institut, Universität Münster, Einsteinstraße 62, 48149 Münster, Germany \\
%              \email{m.kemper@uni-muenster.de}           %  \\
%%             \emph{Present address:} of F. Author  %  if needed
%           J. Lohkamp \at
%              Mathematisches Institut, Universität Münster, Einsteinstraße 62, 48149 Münster, Germany \\
%              \email{j.lohkamp@uni-muenster.de}
%}
\maketitle

\begin{abstract}
Gromov hyperbolic spaces have become an essential concept in geometry, topology and group theory.
Here we extend Ancona's potential theory on Gromov hyperbolic manifolds and graphs of bounded geometry to a large class of Schrödinger operators on Gromov hyperbolic metric measure spaces, unifying these settings in a common framework ready for applications to singular spaces such as RCD spaces or minimal hypersurfaces. Results include boundary Harnack inequalities and a complete classification of positive harmonic functions in terms of the Martin boundary which is identified with the geometric Gromov boundary.
\end{abstract}

\tableofcontents

\phantomsection\addcontentsline{toc}{section}{Introduction}
\section*{Introduction}
Gromov hyperbolic spaces are abstractions of classical manifolds of negative sectional curvature. They are characterized by a large-scale geometric property: all geodesic triangles are uniformly thin. This form of hyperbolicity has a broad relation to group-theoretic or topological problems. In this article, we will focus on a more analytic facet of Gromov hyperbolicity: its striking impact on the potential theory of uniformly elliptic operators. In simple terms it says that the analysis conforms to the large-scale geometry. For example, the analytically defined Martin boundary of such operators coincides with the geometric Gromov boundary. Alano Ancona implemented this appealing theory on manifolds \cite{Anc87,Anc90} and graphs \cite{Anc88} of bounded geometry.

Our aim is to generalize this theory to more singular metric spaces with a less restrictive notion of bounded geometry, unifying amongst others the cases of manifolds and graphs. It turns out that this can be done for complete locally compact separable metric measure spaces with a strongly regular coercive Dirichlet form $\E^0$ as soon as the space is locally doubling and admits local Poincaré inequalities and local lower Dirichlet eigenvalue bounds, with uniform constants (see \autoref{sec:dirichlet} for precise definitions). These mild conditions are satisfied on infinite graphs of bounded degree, complete manifolds with uniform Lipschitz charts onto Euclidean balls, complete manifolds with lower Ricci curvature bounds, and even RCD($K,N$) spaces for $K\in\R$ and $N<\infty$, each with the appropriate version of a Laplacian.

We start by developing the local potential theory on these spaces, i.e., maximum principles, Harnack inequalities, solution of the Dirichlet problem, existence of and local estimates for the Green function. The necessary local estimates for \emph{strongly local} Dirichlet forms, analogous to the Laplacian, were already obtained in this generality by Biroli--Mosco and Sturm in the 1990s \cite{BM95, Stu96}. To include the case of Schrödinger operators, which is important for geometric applications, we further generalize strongly local Dirichlet forms $\E^0$ to \emph{Schrödinger forms} $\E=\E^0+V$ for a not necessarily positive $L^\infty$-function $V$. For coercive Schrödinger forms, most potential theoretic results can be preserved.
This permits us to leverage the power of Brelot's axiomatic potential theory with the tool of balayage and obtain integral representations for positive superharmonic functions. We will recall the relevant concepts in \autoref{sec:pottheo}.

Until now, Gromov hyperbolicity has not entered the picture. This only happens in the final part, where our goal is to identify the abstract potential theoretic Martin boundary with the geometric Gromov boundary. This gives a unique representation of positive harmonic functions as Radon measures on the Martin/Gromov boundary.
With the tool set at hand, Ancona's original arguments can rather swiftly be transferred to the more general case. We try to streamline them and show that constants in a priori estimates depend only on coarse geometric and analytic data such as the constants in the Poincaré inequality and in the doubling, hyperbolicity and coercivity conditions. Our presentation emphasizes the central role of the boundary Harnack inequality on Gromov hyperbolic spaces.

\paragraph{Uniform Spaces}
Gromov hyperbolic geometries and their particular potential theory may appear to be a remote world on their own.
However, in the same vein as the unit disc in $\C$ carries a conformal hyperbolic metric, there are large classes of highly singular spaces which can be conformally deformed to complete Gromov hyperbolic spaces. Now the strategy, also initiated by Ancona \cite{Anc87,Anc90}, is to transfer the potential theory on hyperbolic spaces to spaces admitting such a deformation -- Bonk, Heinonen and Koskela \cite{BHK01} revealed that these are precisely the so-called \emph{uniform spaces}. Special cases of particular interest are \emph{uniform domains} in Euclidean space.

These ideas can be expanded to study even more singular spaces. Our sample case is that of minimal (hyper)surfaces with their extraordinarily difficult singularities. Here we think of the singular set as a boundary of its regular complement. These spaces satisfy a strong form of uniformity, the \si-uniformity, which we use to establish a link to the potential theory of Gromov hyperbolic spaces. This gives us a grip on the otherwise hardly approachable asymptotic analysis towards the singularities from an outsourcing to the Gromov hyperbolic theory. This part is due to work of the second named author \cite{Loh18,Loh19,Loh21}.

\paragraph{Organization of the Paper} In \autoref{sec:dirichlet}, we recall definitions for Dirichlet forms and collect several results on their local potential theory such as maximum principles, Harnack inequalities, solution of the Dirichlet problem, Hölder continuity of harmonic functions and estimates for the Green function.

These results are transferred to Schrödinger forms in \autoref{sec:schrödinger}. For coercive forms, we employ the resolvent equation to derive global exponential decay estimates for the Green function.

The \autoref{sec:pottheo} introduces concepts from axiomatic potential theory such as balayage which we now have at our disposal. The abstract Martin boundary classifies positive harmonic functions, but we have yet to see that it in fact coincides with the Gromov boundary on Gromov hyperbolic spaces.

These are properly introduced in \autoref{sec:gromhyp}, together with the mechanism to generate them from uniform spaces such as uniform Euclidean domains or singular minimal hypersurfaces.

Finally, in \autoref{sec:ident}, we prove the 3G-inequality, which as an intermediate technical result yields the boundary Harnack inequality on Gromov hyperbolic spaces. This in turn easily proves the identification of Martin and Gromov boundary. Transferred to uniform spaces, we can represent harmonic functions as integrals over their (metric) boundary, or over the singular set in the case of minimal hypersurfaces.

\section{Strongly Regular Dirichlet Forms}\label{sec:dirichlet}
Let $(X,d,\mu)$ be a locally compact separable metric measure space and $\E$ a positive semidefinite ($\E(u,u)\ge0$ for all $u\in\F$) symmetric bilinear form on a dense subspace $\F\subset L^2(X,\mu)$ which is a Hilbert space with scalar product $\E_1(u,v):=\E(u,v)+(u,v)$ for $u,v\in\F$, i.e., $\E$ is a \emph{symmetric closed form} on $L^2(X,\mu)$ in the sense of \cite[Definition~I.2.3]{MR92}. Here $(\cdot,\cdot)$ denotes the standard scalar product on $L^2(X,\mu)$.

We call $\E$ \emph{coercive} (on $L^2(X,\mu)$), if there is a $c>0$ such that $\E(u,u)\ge c\cdot(u,u)$ for every $u\in \F$.\footnote{Note that this condition is stronger than $\E$ being a \emph{coercive closed form} in the sense of \cite[Definition~I.2.4]{MR92} (which is not necessarily symmetric) and symmetric, because such a form is only coercive on $(\F,\E_1)$ but not on $(L^2(X,\mu),(\cdot,\cdot))$, see \cite[Remark~I.2.5]{MR92}. We will not use this terminology.}
For such a coercive $\E$, the norms on $\F$ induced by $\E_\alpha(u,v):=\E(u,v)+\alpha\cdot(u,v)$ for $\alpha>-c$ are all equivalent, in particular $\E_\alpha$ is a symmetric closed form on $L^2(X,\mu)$ for any $\alpha>-c$.

For any symmetric closed form $\E$, we get a family of \emph{resolvents} (Green operators) $(G_\alpha)_{\alpha>0}$ which are bounded linear operators defined on $L^2(X,\mu)$ with values in $\F\subset L^2(X,\mu)$ such that
\[
    \E_\alpha(G_\alpha \phi, u)=(\phi,u)\quad\text{for any $\phi\in L^2(X,\mu)$, $u\in\F$.}
\]
Moreover, $(G_\alpha)$ is a \emph{strongly continuous contraction resolvent} (on $L^2(X,\mu)$), i.e., the following conditions hold:
\begin{enumerate}
\item $\lim_{\alpha\to\infty}\alpha G_\alpha u=u$ for all $u\in L^2(X,\mu)$ (strong continuity),
\item $\|\alpha G_\alpha u\|\le \|u\|$ for any $u\in L^2(X,\mu)$ (contraction), and
\item $G_\alpha - G_\beta = (\beta-\alpha) G_\alpha G_\beta$ for $\alpha,\beta>0$ (resolvent equation).
\end{enumerate}

$\E$ is called \emph{regular}, if $\F\cap C_\mathrm{c}(X)$ is dense in $\F$ w.r.t.\ the norm $\sqrt{\E_1}$ and dense in $C_\mathrm{c}(X)$ (continuous functions with compact support) in the topology of uniform convergence.

$\E$ is \emph{local}, if $\E(u,v)=0$ for $u,v\in\F$ with disjoint compact support,
and \emph{strongly local}, if $\E(u,v)=0$ in the case that $u\in\F$ is constant in a neighborhood of the support of $v\in\F$.

\paragraph{Harmonic Functions}
For any open set $U\subset X$, $\E$ is locally defined on $\Floc(U):=\{u\in L^2_\loc(U,\mu)\mid \forall V\Subset U \exists \tilde u\in\F: u|_V\equiv \tilde u|_V \text{ $\mu$-a.e.}\}$. Dually, we can consider the space $\Fc(U)$ of functions in $\F$ with (essentially) compact support in $U$ and its closure $\F_0(U)$ in the norm $\sqrt{\E_1}$. Then for local $\E$, $\E(u,\phi)$ is well-defined for $u\in\Floc(U)$ and $\phi\in \F_0(U)$.

A continuous function $u\in\Floc(U)$ is called ($\E$-)\emph{harmonic} in $U$ if $\E(u,\phi)=0$ for any $\phi\in\Fc(U)$.\footnote{We always take continuous versions of harmonic functions and hence can talk about maxima on compact sets instead of essential maxima etc. We will see later that there are sufficiently many such functions.}

More generally, $u\in\Floc(U)$ is ($\E$-)\emph{superharmonic} in $U$ if it is lower semi-continuous with values in $(-\infty,\infty]$ and $\E(u,\phi)\ge0$ for any $\phi\in\Fc(U)$ with $\phi\ge0$. A function $u$ is \emph{subharmonic} if $-u$ is superharmonic.

Sometimes we will express this in terms of the self-adjoint closed operator $\L$ on $L^2(X,\mu)$ (also called the infinitesimal generator) associated with $\E$ by the relation $\E(u,v)=(\L u,v)$ for $u,v\in\F$ \cite{MR92}:
We will call $\E$-(sub-/super-)harmonic functions $\L$-(sub-/super-)harmonic, and for a function $u\in\Floc(U)$, we write $\L u\ge0$ (on $U$) if it is $\L$-superharmonic, $\L u=0$ (on $U$) if it is $\L$-harmonic, and we will call $u$ \emph{a solution of the equation $\L u=f$} (\emph{on $U$}) for some function $f\in L^2_\loc(U)$ if
\[
    \E(u,\phi)=(f,\phi)=\int_{U} f \phi \,\d\mu\quad\text{for every $\phi\in\F_\c(U)$}.
\]

\paragraph{Dirichlet Forms} A symmetric closed form is called a \emph{Dirichlet form} if it is \emph{Markovian}, i.e., for any $u\in\F$, the function $v=u^+\wedge 1$ is in $\F$, and $\E(v,v)\le\E(u,u)$.

\paragraph{Energy Measure}
By the Beurling--Deny theorem, a strongly local, regular Dirichlet form $\E^0$ on $\F$ can be written as %TODO ref
\[
    \E^0(u,v)=\int_X\d \Gamma(u,v)
\]
where the \emph{energy measure} $\Gamma$ is a positive semidefinite symmetric bilinear form on $\F$ with values in the signed Radon measures on $X$.
$\d\Gamma$ satisfies a Leibniz product rule in both parameters.
Note that we henceforth reserve the notation $\E^0$ for strongly local regular Dirichlet forms to tell them apart from more general symmetric closed forms such as Schrödinger forms considered later. The associated self-adjoint operator will be called $\L^0$, the resolvent $G_\alpha^0$.
\paragraph{Induced Pseudo-Metric and Strong Regularity}
The energy measure of a strongly local, regular Dirichlet form $\E^0$ can be used to define a pseudo-metric $\rho$ on $X$:
\begin{multline*}
    \rho(x,y) := \sup\biggl\{u(x)-u(y)\biggm| u\in \F_\loc(X)\cap C(X),\\\diff{\Gamma(u,u)}{\mu}\text{ exists and is $\le1$ $\mu$-almost everywhere}\biggr\}.
\end{multline*}
The strongly local, regular Dirichlet form $\E^0$ is called \emph{strongly regular}, if this $\rho$ is in fact a metric, the \emph{intrinsic metric}, and induces the same topology as the original metric $d$ on $X$.
As we only used topological properties of $d$ up to this point, and only intend on referring to the intrinsic metric $\rho$, we can (and will) just as well assume $d=\rho$ for simplicity of presentation. A posteriori, one can put another metric on $X$ which is bi-Lipschitz equivalent to $\rho$.
This may be seen as a uniform ellipticity condition, and all results are preserved or require a mere change in constants.

We further require that the metric $d=\rho$ is complete.

For several properties of the intrinsic metric, see \cite{Stu96}. E.g., for a strongly regular Dirichlet form, every closed ball is compact, and if $\rho$ is complete, $(X,\rho)$ is a geodesic metric space, i.e., any two points can be connected by a length-minimizing geodesic \cite[Lemma~1.1]{Stu96}.
%TODO mehr Motivation: energy measure ~ |gradient|²; Bsp. Riemannsche Mfkt.en (Stu96, §3); uniform elliptisch: Metrik d bilipschitz (oder gar quasiisometrisch) zu \rho reicht für die meisten Resultate

\subsection{Bounded Geometry Conditions}\label{sec:bg}
In this section, let $\E^0$ be a strongly regular Dirichlet form on a complete space $X$.

We impose three further conditions that essentially limit the local geometric complexity in a uniform way.
The first property of the metric measure space $X$ locally bounds its (Assouad) dimension from above:
\begin{property}[uniform local doubling]\label{prop:doubling}
There are constants $\sigma>0$ and $N>0$ such that
\[
    \mu(B_{2r}(x))\le 2^N \mu(B_r(x))
\]
for all balls $B_{2r}(x)\subset X$ with $0<r<\sigma$.
\end{property}

The other two conditions ensure that functions can be uniformly controlled by their derivatives as represented by the energy measure $\Gamma$ associated with $\E^0$.
\begin{property}[uniform local Poincaré inequality]\label{prop:poincare}
There are constants $\sigma>0$ and $C_P>0$ such that
\[
    \int_{B_r(x)}|u-u_{B_r(x)}|^2\,\d\mu\le C_P r^2\int_{B_{2r}(x)}\d\Gamma(u,u)\text{ for $u\in \Floc(X)$}
\]
for all balls $B_{2r}(x)\subset X$ with radius $0<r<\sigma$.
\end{property}

\begin{property}[uniform local Dirichlet eigenvalue bound]\label{prop:dirichlet}
There are constants $\sigma>0$ and $C_D>0$ such that
\[
    \int_{B_r(x)}|u|^2\,\d\mu\le C_D r^2\E^0(u,u)\text{ for $u\in \F_0(B_r(x))$}
\]
for all balls $B_{r}(x)\subset X$ with radius $0<r<\sigma$.
\end{property}

Since $X$ is assumed complete, the doubling property \ref{prop:doubling} implies \emph{upper} bounds $\sim r^{-2}$ for the lowest non-zero eigenvalues of the self-adjoint operator $\L^0$ associated with $\E^0$ on balls of radius $r<\sigma$ with Neumann as well as Dirichlet boundary conditions, while the other two properties yield \emph{lower} bounds with the same asymptotic for Neumann (\autoref{prop:poincare}) and Dirichlet (\autoref{prop:dirichlet}) eigenvalues, respectively \cite[Section~2.A)]{Stu96}. Although \cite{Stu96} does not need the lower Dirichlet eigenvalue bounds (\enquote{local spectral gap}) \ref{prop:dirichlet} to derive Harnack inequalities, we want to apply stronger quantitative results from \cite{BM95} relying on these bounds, cf.\ the discussion in \cite[Section~2.B)]{Stu96}.

At the first glance, these conditions seem to fix a scale $\sigma$ up to which they hold, but this scale can in fact be chosen arbitrarily large (with a subsequent loss in the other constants; see e.g.\ \cite[Lemma 4.11]{BB18} for Poincaré inequalities. For the doubling condition, it is essential that $X$ is geodesic). This sets our bounded geometry conditions apart from more restrictive ones such as lower bounds on the injectivity radius of the exponential maps on Riemannian manifolds (plus curvature bounds) or uniformly bi-Lipschitz maps to equally sized balls in Euclidean space.

The conditions are satisfied in the following examples:
\begin{itemize}
  \item In Euclidean space $\R^n$ with Laplacian $\L^0=-\Delta$, they are known to hold for $\sigma=\infty$.
  \item They are satisfied on finite metric graphs with the graph Laplacian. For infinite graphs, they e.g.\ follow from uniform bounds on edge lengths and the degree of vertices.%TODO ref
  \item On closed Riemannian manifolds with the Laplacian, they are true by local comparison with Euclidean space.
  \item On complete noncompact Riemannian manifolds, one needs additional bounded geometry conditions. Uniform bounds on sectional curvature and injectivity radius work, but in fact a lower bound on the Ricci curvature is sufficient.
  \item More general, one can consider $\mathrm{RCD}^*(K,N)$ spaces for finite $N$ and arbitrary $K\in\R$ with their Cheeger energy $\mathrm{Ch}_2$, see \cite{Gig15} and the references therein. %TODO better ref; Hajjlasz--Koskela for (1,2)->(2,2) Poincaré? -> \cite[Theorem~4.18]{Hei01}
      The Dirichlet eigenvalue bound can be seen as follows: Under blow-ups, subsequences converge to Euclidean space \cite{GMR15}, where such a (scaling-invariant) bound is known to hold. As Dirichlet eigenvalues on balls are lower semi-continuous (see e.g.\ \cite{AH18}), a series of counterexamples to the eigenvalue bound would produce a contradiction.
  \item In all these examples, the metric induced by the Dirichlet form agrees with the original metric, but one could easily consider other \emph{uniformly elliptic} operators than the Laplacian with a merely bi-Lipschitz equivalent induced metric. In the end, qualitative results will carry over while quantitative results in terms of the metric can easily be adjusted.
\end{itemize}

\subsection{Local Bounds and Dirichlet Problem}
\label{sec:strreglocbounds}
For a strongly regular Dirichlet form $\E^0$ with associated self-adjoint operator $\L^0$ on a complete doubling space with Poincaré inequalities and Dirichlet eigenvalue bounds, there are explicit bounds on heat kernels and Green functions in terms of the constants $\sigma$, $N$ and $C_P$ \cite{BM95,Stu96}. They can be used to easily prove Harnack inequalities, Hölder continuity of harmonic functions and to solve the Dirichlet problem. Here, we state the more general bounds which will later enable us to obtain these corollaries as well in the case of Schrödinger forms.

\paragraph{Minimum Principle}
The minimum principle asserts that if $U\subset X$ is open and bounded and $u$ a $\L^0$-superharmonic function on $U$ extending lower semi-continuously to $\del U$, then $u\ge \min_{\del U}u$ on $U$ \cite[2.(i)]{BM95}. In particular, if $u$ is nonnegative on the boundary, it is nonnegative everywhere.

Analogously, there is a maximum principle for $\L^0$-subharmonic functions. $\L^0$-harmonic functions satisfy both these principles.

\paragraph{Dirichlet Problem with Homogeneous Boundary Condition}
The problem $\L^0u=\phi$ on a ball $B$ with homogeneous boundary condition $u\in\F_0(B)$ has well-behaved solutions for suitable functions $\phi$. For later use, we state the following crucial estimate:
\begin{theorem}[$L^\infty$-estimate on balls]\label{thm:Linf}
\cite[Theorem~4.1]{BM95}
Let $x\in X$ be a point, $0<2r<\sigma$ and $u\in\F_0(B_r(x))$ a solution of the equation $\L^0 u = \phi$ on $B_r(x)$

for some $\phi\in L^p(B_r(x))$, $p>\max(N/2, 2)$. Then there is a constant $c_\infty=c_\infty(N, C_P, C_D)$ such that
\[
\|u\|_{L^\infty(B_r(x))}\le  \frac12 c_\infty r^2(\mu(B_r(x))^{-1/p} \|\phi\|_{L^p(B_r(x))}\le \frac12 c_\infty r^2 \|\phi\|_{L^\infty(B_r(x))}.
\]
\end{theorem}

In the situation of this Theorem, solutions exist and are unique by the Poincaré inequality.
Their regularity follows from a more general result:

\paragraph{Regularity of Local Solutions}
On any open set $U\subset X$, solutions $u\in \F_\loc(U)$ of the equation $\L^0u=\phi$ for some $\phi\in L^p(B_r(x))$, $p>\max(N/2, 2)$, are locally Hölder continuous with constants depending only on $\sigma$, $N$, $C_P$ and $C_D$ \cite[Theorem~5.13]{BM95}. This is in particular true for $\L^0$-harmonic functions.

\paragraph{Local Green Functions}
For every ball $B_r(x)$ with $0<2r<\sigma$ and point $y\in B_r(x)$, there is a unique positive Green function $G^0_{B_r(x)}(\cdot, y)\in L^{p'}(B_r(x),\mu)$ which satisfies
\[
    \int_{B_r(x)}G^0_{B_r(x)}(\cdot, y)\,\phi\,\d\mu=u(y),
\]
for any solution $u$ of $\L^0u=\phi$ as in \autoref{thm:Linf}, see \cite[6]{BM95}. It is $\L^0$-superharmonic on $B_r(x)$ and $\L^0$-harmonic on $B_r(x)\setminus\{y\}$, and hence locally Hölder continuous on $B_r(x)\setminus\{y\}$ by the preceding paragraph. Note that for general balls, the Green function does not necessarily vanish towards the boundary. %TODO example?

\begin{theorem}[estimate for Green functions on balls]\label{thm:g0locest}
(\cite[Theorem~1.3]{BM95}; see also \cite[Cor.~4.11]{Stu96})
The Green function $G^0_{B_{r}(x)}(\cdot, x)$ of $\E^0$ on a ball $B_{r}(x)$, for any $x\in X$ and $0<20r<\sigma$, satisfies
\[
C_G^{-1}\int_{d(x,y)}^r\frac{s^2}{\mu(B_s(x))}\frac{\d s}s \le G^0_{B_{r}(x)}(y,x)\le C_G \int_{d(x,y)}^r\frac{s^2}{\mu(B_s(x))}\frac{\d s}s \quad\text{for any $y\in B_{r/16}(x)$}
\]
with a constant $C_G\ge 1$ only depending on $N$, $C_P$ and $C_D$.
\end{theorem}

\paragraph{Dirichlet Problem with Inhomogeneous Boundary Condition}
A more general boundary value problem is, given an open set $V\Subset X$ with $\del V\neq\emptyset$, to find a function $u\in C(\overline V)$ which is $\L^0$-harmonic on $V$ and satisfies $u|_{\del V}=f$ for a prescribed function $f\in C(\del V)$.
Such a set $V$ is called \emph{regular} if the Dirichlet problem admits a solution for any continuous boundary value.

Solutions of the Dirichlet problem are always unique by the minimum principle.

In general, balls do not have to be regular, but nonempty sets of the form
\[
    U_{-\delta}:=\{x\in U\mid \dist(x,X\setminus U)>\delta\}
\]
with $\delta>0$ and $U\Subset X$, $\del U\neq\emptyset$, are always regular: %TODO example balls, S^n
Since $X$ is geodesic, these sets have an exterior corkscrew and hence satisfy the Wiener criterion \cite[Lemma~14.2]{BB11}.
The Wiener criterion for solvability of the Dirichlet problem for $\E^0$ is valid in our situation by \cite{Bir94,BM07}, hence we can apply the construction above to balls to obtain a topological basis of regular sets as well as a compact exhaustion by regular sets in the case that $X$ is unbounded \cite[Theorem~14.1]{BB11}.

On regular sets, we can get the following version of \thmref{thm:Linf}:
\begin{corollary}[$L^\infty$-estimate on regular sets]\label{thm:Linfreg}
Let $U\Subset X$ be a regular set with $\diam(U)<\sigma/2$ and $\phi\in L^\infty(U)$. Then there is a unique solution $u\in C(\overline U)$ of the equation $\L^0u=\phi$ with boundary condition $u|_{\del U}\equiv 0$ which satisfies
\[
\|u\|_{L^\infty(U)} \le c_\infty (\diam U)^2 \|\phi\|_{L^\infty(U)}
\]
with constant $c_\infty=c_\infty(N, C_P, C_D)$ from \autoref{thm:Linf}.
\begin{proof}
First we solve the equation $\L^0 u'=\phi$ on a ball $B$ of radius $\diam(U)$ containing $U$, where $\phi$ is extended by zero outside $U$. Then we can subtract the Dirichlet solution on $U$ with boundary value $u'|_{\del U}$ to obtain the desired $u$. The bounds on $u'$ from \thmref{thm:Linf} carry over by the maximum principle with an additional factor 2.
\end{proof}
\end{corollary}

\paragraph{Harnack Inequality}
Finally, in our setting there is an elliptic Harnack inequality for nonnegative harmonic functions.
\begin{theorem}[elliptic Harnack inequality]\label{thm:harnack0}\cite[Proposition~3.2]{Stu96}, \cite[Theorem~1.1]{BM95}
There is a constant $H\ge 1$ only depending on $\sigma$, $N$ and $C_P$ such that any \emph{nonnegative} solution $u$ of $\L^0 u=0$ on a ball $B_{2r}(x)$ with $0<r<\sigma$ satisfies
\[
    \sup_{B_r(x)} u \le H\cdot \inf_{B_r(x)} u.
\]
\end{theorem}
This directly implies a version of the strong minimum principle: If a nonnegative $\L^0$-harmonic function on a connected open set takes the value 0 in some point, then it is $\equiv 0$.

\section{Schrödinger Forms}\label{sec:schrödinger}
Henceforth, we assume that $\E^0$ is a strongly regular Dirichlet form on a complete space $X$ satisfying the bounded geometry conditions from \autoref{sec:bg}.

For any not necessarily positive function $V\in L^\infty(X,\mu)$, %\emph{signed} Radon measure $\nu$ on $X$,
the regular local symmetric closed form
\[
    \E(u,v):=\E^0(u,v)+\int_Xuv\,V\,\d\mu
\]
is defined for $u,v\in\F$ with the same domain as $\E^0$. We will call this a \emph{Schrödinger form}, since its associated self-adjoint operator $\L=\L^0+V$ is a Schrödinger operator.
Note that $\E$ might be neither strongly local (for any nonzero $V$) nor Markovian (if $V$ is somewhere negative).

This extends our set of global constants by an essential upper bound for $V$, we will call this $k\ge 0$ and assume
\[
    \|V\|_\infty\le k.
\]
Examples from quantum mechanics, where Schrödinger operators originated, suggest to study potentials $V$ with singularities, e.g., in (local) $L^p$ spaces. Instead, we opted for the locally most simple case $L^\infty$ to not further complicate the already technical local analysis in \autoref{sec:locbounds}. Starting from \autoref{sec:weakcoercivity}, there will be no more explicit reference to $V$, hence any generalization giving the same local results (solution to Dirichlet problem, bounds on Green function, Harnack inequality) is also a viable starting point for the subsequent global analysis.

\subsection{Local Bounds}\label{sec:locbounds}
In this section, we want to transfer the bounds from \autoref{sec:strreglocbounds} and their standard consequences (Harnack inequalities, Hölder continuity of harmonic functions) to Schrödinger forms $\E=\E^0+V$.

The results in this section a priori hold only on small balls with a maximal radius depending on $\sigma$, $N$, $C_P$, $C_D$ and $k$. To avoid the clutter of an additional constant, we simply assume that they are already satisfied on balls of radius $\sigma$, which can always be achieved by choosing a smaller value of $\sigma$. With the additional coercivity assumption in the next section in place, the results in this section actually hold for arbitrarily large balls, so that for later applications, we recover these results with the original $\sigma$ anyway.

\subsubsection{Dirichlet problem}
We start with a solution $u_0$ of $\L^0u_0=\phi$ for some $\phi\in L^\infty(U)$ on a regular set $U$ with boundary condition $u_0|_{\del U}\equiv0$.
Now we iterate as follows: We let $u_i$ be the solution of $\L^0\,u_i=-V\,u_{i-1}$ with boundary condition $u_i|_{\del U}\equiv 0$. Then $\L(u_0+\cdots+u_i)=\phi+V\,u_i$. By the \thmnameref{thm:Linfreg}, there is a constant $c_\infty=c_\infty(\sigma, N, C_P, C_D, k)>0$ such that
\[
    \sup_{U}|u_i|\le c_\infty (\diam U)^2 \|V\,u_{i-1}\|_{L^\infty(U)},
\]
and hence for $(\diam U)^2<(c_\infty \|V\|_{L^\infty(X)})^{-1}$, the series $\sum_{i=0}^\infty u_i$ converges uniformly on $\overline{U}$ to a limit function $u$ with $\L\,u=\phi$ and $u|_{\del U}\equiv0$. %TODO convergent in which sense to preserve operator?
Moreover, for $(\diam U)^2<\frac12(c_\infty \|V\|_{L^\infty(X)})^{-1}$, we have the estimate
\[%begin{equation}\label{eq:LinfSchrödinger}
    \|u\|_{L^\infty(U)} \le \|\phi\|_{L^\infty(U)} \sum_{i=1}^\infty (c_\infty (\diam U)^2 \|V\|_{L^\infty(X)})^i
    \le 2c_\infty (\diam U)^2 \|V\|_{L^\infty(X)} \|\phi\|_{L^\infty(U)}.
\]%end{equation}

Now we can easily solve Dirichlet problems with arbitrary continuous boundary conditions on such a small regular set $U$: Given a continuous function $f$ on $\del U$, there is a $u_0$ with $\L^0u_0=0$ and $u_0|_{\del U}=f$ and a $v$ with $\L\,v=-V\,u_0$ and $v|_{\del U}\equiv 0$, such that $u:=u_0+v$ solves $\L\,u=0$ and $u|_{\del U}= f$.

Solutions to all these problems are \emph{unique} on sufficiently small balls. This follows from the \thmnameref{prop:dirichlet} applied to the difference $u\in\F_0(B_r(x))$ of two solutions: for $r^2\le1/(C_D\|V\|_\infty)$, this yields
\begin{align*}
    0=\E(u,u)&=\E^0(u,u)+\int_{B_r(x)}u^2\,V\,\d\mu\ge\E^0(u,u)-\|V\|_\infty\int_{B_r(x)}u^2\,\d\mu\\
    &\ge \left(\frac1{C_Dr^2}-\|V\|_\infty\right)\int_{B_r(x)}u^2\d\mu\ge 0
\end{align*}
and hence $u\equiv0$.

\subsubsection{Local Green Functions}\label{sec:locgreen}
We note that the measure $\mu$ has a strictly positive \emph{lower regularity dimension}: there are constants $C_L, \alpha>0$ such that
\[
\mu(B_r(x))\ge C_L\left(\frac{r}{\rho}\right)^\alpha \mu(B_\rho(x))\quad\text{for any $x\in X$ and $0<\rho<r<\sigma$.}
\]
This follows from the doubling property and the fact that $X$ is connected and unbounded, see e.g.\ \cite[Proposition~5.2]{GH14}, and the constants only depend on $\sigma$ and the doubling constant $N$. We may assume $0<\alpha<1$.

Using this, the \thmnameref{thm:g0locest} gives for a Green function $G^0_{B}$ on a ball $B_r:=B_r(x_0)$, $0<20r<\sigma$, and $0<2\rho\le r/16$:
\begin{multline*}
  \int_{B_{2\rho}\setminus B_\rho}G^0_{B}(y,x_0)\,\d\mu(y)\mule \mu(B_{2\rho})\int_\rho^r\frac{s\,\d s}{\mu(B_s)}
  \mule \frac{\mu(B_{2\rho})}{\mu(B_{\rho})}\int_\rho^r\frac{s\,\d s}{(s/\rho)^\alpha}\\
  \mule\rho^\alpha(r^{2-\alpha}-\rho^{2-\alpha})\le \rho^\alpha r^{2-\alpha},
\end{multline*}
where we use the notation \enquote{$\mule$} meaning \enquote{up to a constant factor depending only on $\sigma$, $N$, $C_P$, $C_D$ and $k$}.
Hence,
\[
    \int_{B_{r/16}}G^0_{B}(y,x_0)\,\d\mu(y) \mule\sum_{i=0}^\infty \frac{r^2}{2^{\alpha i}}\mule r^2.
\]
Moving the pole and comparing Green functions on different balls, we can get the more flexible version
\[
    \int_{B_{r/32}}G^0_{B'}(y,x)\,\d\mu(y)\mule r^2\quad\text{for any $x\in B_{r/32}(x_0)$}
\]
on $B':=B_{r/2}(x_0)$.

We will use this bound to construct a local Green function for $\E$ with an iteration similar to the Dirichlet problem:

As distributions, we have $\L_xG^0_{B'}(x,y)=\delta_y(x)+W(x,y)$, with $W(x,y)=V(x)G^0_{B'}(x,y)$. By the above bounds and global bounds on $V$, we can find a small radius $r_{\mathrm{gb}}$ only depending on $\sigma$, $N$, $C_P$, $C_D$ and $k$ such that $\int_{B_{r_\mathrm{gb}}}|W(x,y)|\,\dy<1/2$ for $x\in B_{r_\mathrm{gb}}$. This implies the finiteness of all iterates $W^i$ where the product of kernels $S$, $T$ is defined as $(ST)(x,y)=\int_{B_{r_\mathrm{gb}}}S(x,z)\,T(z,y)\,\d\mu(z)$. The bound ensures also the summability of
\[
    \tilde G:=G^0_{B'}\sum_{i=0}^\infty(-W)^i\,,
\]
again as products of kernels.
Application of $\L$ yields
\[
    \L\tilde G=(\delta+W)\sum_{i=0}^\infty(-W)^i=\sum_{i=0}^\infty(-W)^i-\sum_{i=1}^\infty(-W)^i=\delta\,,
\]
so that $\tilde G$ is indeed a Green function for $\L$ on $B_{r_\mathrm{gb}}(x_0)$. The bounds on $W$ directly give local estimates on $\tilde G$ as in \thmref{thm:g0locest}.

\subsubsection{\texorpdfstring{\textit{h}}{h}-Transform and Harnack Inequality}
We will utilize Doob's $h$-transform technique as a means to generalize results for strongly local Dirichlet forms to Schrödinger forms, see \cite{GS11} for more details.

\begin{lemma}[$h$-transform]\label{lem:htrafo}
For a domain $U\subset X$ and a continuous function $h>0$, the \emph{$h$-transform} of a Schrödinger form $\E$ is the local regular symmetric closed form $\E^h(u,v):=\E(u\cdot h, v\cdot h)$ on $\F^h:=\F_0(U)/h\subset L^2(U,h^2\mu)$.

It has the following properties:
\begin{enumerate}
    \item If $\E$ is coercive on $U$, i.e., there is a $c>0$ such that $\E(u,u)>c\cdot(u,u)$ for all $u\in\F_0(U)$, then $\E^h$ is coercive with respect to the scalar product on $L^2(U,h^2\mu)$ with the same $c$.
    \item If $h$ is $\E$-superharmonic on $U$, then $\E^h$ is a Dirichlet form.
    \item If $h$ is $\E$-harmonic and bounded from above and away from 0 on $U$, $\E^h$ induces the same metric on $U$ as $\E^0$, and $\E^h$ satisfies properties \ref{prop:doubling}--\ref{prop:dirichlet} on balls of radius $<\sigma$ contained in $U$.
\end{enumerate}
\begin{proof}
    \begin{enumerate}
        \item Follows directly from the definition.
        \item For any $u'=u/h\in\F^h$, we have $u'^+\wedge 1 = (u^+\wedge h)/h$ and hence
        \[
            \E^h(u'^+\wedge 1, u'^+\wedge 1) = \E(u^+\wedge h, u^+\wedge h) \le \E(u,u)=\E^h(u',u')
        \]
        since $h$ is $\E$-superharmonic.
        \item The bounds on $h$ ensure that $\F_c(U)\cap L^\infty(U,\mu)$ is dense in both Hilbert spaces $\F_0(U)$ and $\F^h$ \cite[Proposition~5.7]{GS11}, hence it is sufficient to consider these functions.
        Using the Leibniz rule for $\d\Gamma$, we have for any $u\in\F_c(U)\cap L^\infty(U,\mu)$
        \begin{align*}
            \E^h(u,u)&=\int_U\left(\d\Gamma(hu,hu)+h^2u^2V\d\mu\right)\\
            &=\int_U\left(\d\Gamma(h,hu^2)+h\cdot hu^2\,\d\mu+h^2\d\Gamma(u,u)\right)=\int_Uh^2\d\Gamma(u,u),
        \end{align*}
        where the last step used that $h$ is $\E$-harmonic and $hu^2\in\F_c(U)$.
        This shows that $h^2\Gamma$ is the energy measure of the strongly local Dirichlet form $\E^h$. In particular, with respect to the measure $h^2\mu$, it induces the same metric on $U$ as $\E^0$.

        On balls of radius $<\sigma$ contained in $U$, the doubling property \ref{prop:doubling} holds with a constant modified by the bounds on $h$. The same is true for the Poincaré inequality \ref{prop:poincare} and the Dirichlet eigenvalue bound \ref{prop:dirichlet} by the expression for the energy measure above, since the involved function spaces share a common dense set.
    \end{enumerate}
\end{proof}
\end{lemma}

This helps us to transfer the Harnack inequality to $\L$-harmonic functions.

\begin{theorem}[Harnack inequality]
\label{harnackineq}
There is a constant $H=H(N, C_P, C_D, k)\ge 1$ such that if $u\ge0$ is $\L$-harmonic on $B_{2r}(x_0)$, for some $0<r\le \sigma$ and $x_0\in X$, then on $B:=B_r(x_0)$,
\[
    \sup_B u \le H \inf_B u.
\]
\begin{proof}
    We will prove the Harnack inequality first for balls $B$ of sufficiently small radius $r$ and then generalize to radii up to $\sigma$.

    In the case $r<r_\mathrm{gb}/3$, there is a harmonic function $h>0$ on $U:=B_{2r}(x_0)$ -- take the Green function on $B_{r_\mathrm{gb}}(x_0)$ with pole in a controlled distance from $U$ as constructed in \autoref{sec:locgreen}. Moreover, if properly scaled with $\mu(B_{r_\mathrm{gb}}(x_0))$, the bounds obtained there universally bound $h$ on $U$ from above and away from 0. Hence, we can apply the $h$-transform from \autoref{lem:htrafo} to get a strongly regular Dirichlet form $\E^h$ satisfying locally the requirements for the \thmnameref{thm:harnack0}. This applies to the $\E^h$-harmonic function $u/h$, giving us a Harnack inequality for $u$ with constants depending only on the Harnack constant from \autoref{thm:harnack0} and the universal bounds on $h$.

    For larger radii, we can apply the small radius version along a \emph{Harnack chain} of controlled length, i.e., a sequence of overlapping small balls where the Harnack inequality applies along geodesics connecting two points with the center $x_0$. Since the number and size of these balls depends only on universal constants, the same is true for the resulting Harnack constant $H$.
\end{proof}

\end{theorem}

A standard consequence is that $\L$-harmonic functions are Hölder continuous, with local estimates for their oscillation in terms of universal constants. As before, the Harnack inequality implies a strong minimum principle:
If a nonnegative $\L$-harmonic function on a connected open set takes the value 0 in some point, then it is $\equiv 0$.

\subsection{Coercivity}\label{sec:weakcoercivity}
If a Schrödinger form $\E$ is \emph{coercive}, the resolvent $G_\alpha$ exists even for $\alpha=0$, and we simply write $G:=G_0$. We still have to see that this operator has an integral kernel, i.e., a global Green function also called $G$ such that
\[
    G(f)=\int_XG(\cdot,y)\,f(y)\,\d\mu(y)\quad\text{for $f\in L^2(X,\mu)$, $\mu$-a.e.}
\]

We start with a version of the weak minimum principle, which utilizes a strictly positive superharmonic function in the $h$-transform technique. Locally, such functions exist, as already used in the proof of the \thmnameref{harnackineq}.
Bootstrapping the resulting local version of the minimum principle, we will show the existence of a global positive superharmonic function for coercive Schrödinger operators in the next step, establishing the minimum principle for arbitrarily large domains.

\begin{proposition}[weak minimum-zero principle]\label{prop:min0}
Let $\E$ be a Schrödinger form which is coercive on $U\Subset X$ and assume there is a continuous $\E$-superharmonic function $h>0$ on a neighborhood of $U$. Let $u$ be an $\E$-superharmonic function on $U$ with $u|_{\del U}\ge 0$. Then $u\ge 0$ on all of $U$.
\begin{proof}
By \thmref{lem:htrafo}, $\E^h$ is coercive on $U$ and (in contrast to $\E$) a Dirichlet form, since $h$ is $\E$-superharmonic.
This means we can apply the weak minimum principle as given in \cite[Proposition~4.6]{GH08} to $\E^h$, yielding the assertion for $\E$ as well.
\end{proof}
\end{proposition}

\begin{lemma}
For any \emph{coercive} Schrödinger form $\E$, there is a continuous positive $\E$-su\-per\-har\-monic function $v$ on $X$.
\begin{proof}
First we note that there is a strictly positive continuous function $\phi\in\F$. To construct such a function, we  take any ball $B$ and a compactly supported continuous function which is 1 on $B$. By regularity, there is a function in $\F\cap C_\c(X)$ which $1/2$-approximates this function, hence is at least $1/2$ on $B$, and can be assumed positive everywhere (using that $\E^0$ is Markovian). Such functions associated to balls in a locally finite countable cover of $X$ can be summed with proper weights to ensure convergence in $\F$, resulting in a strictly positive function $\phi\in\F\cap C(X)$.

Coercivity means in particular that $\F$ is a Hilbert space with scalar product $\E$. By the Riesz representation theorem applied to the continuous functional $\psi\mapsto(\phi,\psi)$ with a fixed $\phi\in \F$ with $\phi>0$ as above, we obtain a $v\in\F$ with
\[
    \E(v,\psi)=(\phi,\psi)\quad\text{for all $\psi\in\F$,}
\]
i.e., $\L v=\phi>0$ for the associated operator. Moreover, as a solution of this equation, $v$ is continuous.
$v$ is nonnegative because for $v^-:=\min\{0,v\}$,
\[
    0\le \E(v^-,v^-)=\E(v,v^-)=\int_X\phi\,v^-\,\d\mu\le0
\]
by locality of $\E$.

To show that $v$ is strictly positive, we want to apply the strong minimum principle stated after the \thmnameref{harnackineq}. This holds for nonnegative harmonic functions, but extends to continuous superharmonic functions as follows: on a small regular ball, the $\E$-harmonic function $u$ with boundary value given by $v$ satisfies $u\le v$ (this is the axiomatic characterization of superharmonic functions we will revisit in \autoref{sec:pottheo}).
Moreover, on a sufficiently small ball, there is a strictly positive continuous superharmonic function $h$, the Green function on a slightly larger ball with pole outside the smaller ball as constructed in \autoref{sec:locgreen}. Hence, we can apply the \thmnameref{prop:min0} to see that $u\ge 0$ on the small ball. Now the strong minimum principle after the \thmnameref{harnackineq} shows that either $u>0$ and hence $v>0$ in the interior of the small ball, or $u\equiv 0$. The latter case cannot happen on all balls with a prescribed center, as this would mean $v\equiv0$ on an open set, which is excluded by the equation $\L v = \phi >0$. Consequently, $v>0$ everywhere.
%Since $v$ satisfies $\L v=\phi>0$, there is a perturbation $h$ of $v$ that is strictly positive and still $\E$-superharmonic and continuous: We can add small, compactly supported nonnegative bump functions to $v$ to make it positive.
\end{proof}
\end{lemma}

As mentioned before, this allows us in the case of a coercive Schrödinger operator to apply \thmref{prop:min0} on any $U\Subset X$. In particular, in this case solutions to Dirichlet problems on relatively compact sets are always unique, since the minimum principle shows that their difference vanishes.

\begin{theorem}[existence of a Green function]\label{thm:green}
If $\E$ is a Schrödinger form which is coercive on a domain $\Omega\subset X$, i.e., the principal Dirichlet eigenvalue
\[
    \lambda_1(\E, \Omega):=\inf_{u\in\F_0(\Omega)}\frac{\E(u,u)}{(u,u)}
\]
is strictly positive (recall that $\F_0(\Omega)$ denotes the closure of the space of functions in $\F$ compactly supported in $\Omega$), then there is a \emph{Green function} for $\E$, i.e., a function $G^\Omega:\Omega\times\Omega\to(0,\infty]$ such that
\begin{enumerate}
\item $G^\Omega$ is harmonic, and in particular continuous, in both arguments away from the diagonal,
\item $G^\Omega(x,y)=G^\Omega(y,x)$, and
\item $(G^\Omega f)(x)=\int_\Omega G^\Omega(x,y)f(y)\,\d\mu(y)$ for any $f\in L^2(\Omega,\mu)$ and a.e.\ $x\in\Omega$.
\end{enumerate}
\begin{proof}
The proof of the similar result \cite[Lemma~5.2]{GH14} applies in our situation. The differences in the requirements are that \cite{GH14} is concerned with coercive, strongly local, regular Dirichlet forms on \emph{globally} doubling spaces, i.e., the doubling property holds for arbitrarily large balls. The latter is not an issue since doubling is only ever applied on balls smaller than a certain radius.
That $\E$ is assumed Markovian is not explicitly used in the proof of \cite[Lemma~5.2]{GH14}, it is sufficient to have uniform local Harnack inequalities and Hölder estimates. %TODO ref
\end{proof}
\end{theorem}

\subsection{Global Bounds and Relative Maximum Principle}\label{sec:glob}
The basic idea of this section is to combine the local estimates we derived in the previous sections with the resolvent equation for Green functions viewed as operators.

The resolvent equation from \autoref{sec:dirichlet} will be used in the following form:
\begin{lemma}[resolvent equation]
\label{resolvent} For any $\lambda>0$ such that $\E_\lambda$ is still coercive, the resolvent $G_{-\lambda}$ as a bounded operator on $L^2(X,\mu)$
satisfies
\[G_{-\lambda}=G+\lambda \cdot G\circ G_{-\lambda}\,.\]
Since all involved operators are positive, this yields the inequalities
\begin{align}
  G & \le G_{-\lambda} \quad\text{and}\label{eq:resineq1}\\
  G\circ G_{-\lambda} & \le \frac{1}{\lambda}G_{-\lambda}\,.\label{eq:resineq2}
\end{align}
In particular, these results hold for the operators acting on characteristic functions of bounded measurable sets (which lie in $L^2(X,\mu)$) and for their kernels, the Green functions, where the resolvent equation has the form
\[
    G_{-\lambda}(x,z)=G(x,z)+\lambda \int_XG(x,y)\,G_{-\lambda}(y,z)\,\d\mu(y).
\]
This follows directly from the characterization of the Green function in \autoref{thm:green} (iii) and Fubini's theorem.
\end{lemma}

We fix an $\epsilon>0$ such that the Schrödinger form $\E_{-\epsilon}$ is still coercive.

\subsubsection{Behavior of Green Functions}
We combine the resolvent equation with Harnack inequalities to derive growth estimates for Green functions. This is completely analogous to \cite[Proposition~4.1]{Anc90}.
\begin{proposition}[bound for the Green function]
\label{prop:bgreen}
There is a universal constant $c_1(\sigma,N,C_P,C_D,k,\epsilon)\ge 1$ such that for the Green function $G$ of $\E$ on $X$,
\begin{align*}
c_1^{-1}/\mu(B_\sigma(y))\le\ &G(x,y) &&\text{if }d(x,y)\le \sigma\text{, and}\\
&G(x,y)\le c_1/\mu(B_\sigma(y)) &&\text{if }d(x,y)=\sigma\,.
\end{align*}
\begin{proof}
The lower bound is directly obtained from the bounds in the end of \autoref{sec:locgreen}, because we have $G(\cdot,y)\ge g(\cdot,y)$ for any Green function $g$ on a smaller domain, and iterated application of the Harnack inequality along a Harnack chain as in the proof of the \thmnameref{harnackineq}.

\begin{figure}
    \centering
    \includegraphics{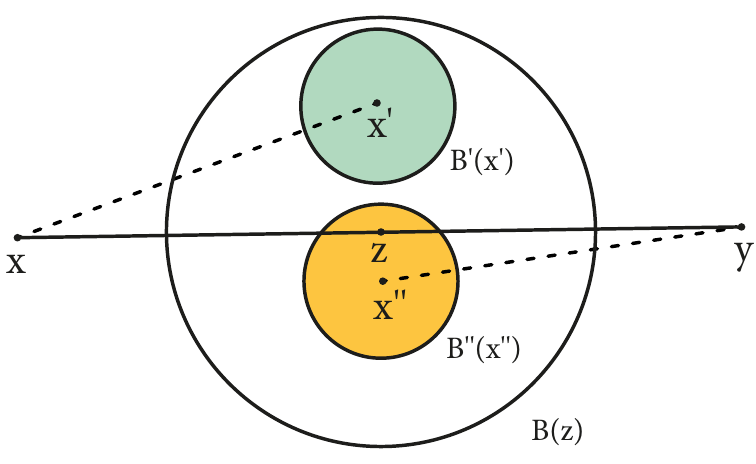}
    \caption{Linking $x$ and $y$ via $B'(x')$ and $B''(x'')$, balls of radius $\sigma/9$, both contained in the ball $B(z)$ of radius $\sigma/3$.
    Along the dashed lines we use Harnack chains to get bounds for the Green function in the proof of \thmref{prop:bgreen}.}
    \label{fig:bgreen}
\end{figure}

For the upper bound, the strategy is to first obtain \emph{any} localized bound and then transfer it to two given points using Harnack chains. To this end, we apply the \thmnameref{resolvent} to the characteristic function $\chi_B$ of an arbitrary ball $B$ of radius $\sigma/3$ (see \autoref{fig:bgreen}) and integrate the result over $B$ to see
\begin{align}\label{eq:respos}
    0 \le\int_X\chi_BG(\chi_B)\,\d\mu &= \int_X\chi_B G_{-\epsilon}(\chi_B)\,\d\mu-\epsilon\int_X\chi_BG(G_{-\epsilon}(\chi_B))\,\d\mu\notag\\ &= \int_X\left[\chi_B-\epsilon G(\chi_B)\right]G_{-\epsilon}(\chi_B)\,\d\mu
\end{align}
using Fubini's theorem. Recall that the Green function $G$ is self-adjoint on $L^2(X,\mu)$.

From \eqref{eq:respos}, we infer the existence of an $x'\in B$ with $\left(\chi_B-\epsilon G(\chi_B)\right)(x')\ge 0$ and therefore
\[ \int_BG(\zeta,x')\,\d\mu(\zeta)=G(\chi_B)(x') \le 1/\epsilon\,. \]
Using that $X$ is geodesic and contains a point outside $B$, we can find another ball $B'\subset B$ of radius $\sigma/18$ (centered somewhere on a geodesic connecting the center of $B$ to the outside), such that $B'\cap B_{\sigma/18}(x')=\ep$. Obviously $\int_{B'}G(\zeta,x')\d\zeta\le\int_{B}G(\zeta,x')\d\zeta\le 1/\epsilon$ and hence there is an $x''\in B'$ with
\[
    G(x'',x') \le 1/(\epsilon\mu(B'))\,.
\]
The result is an upper bound on $G(x'',x')$ for two points with distance at least $\sigma/9$ in an arbitrary ball $B$.

To get an upper bound for any two points $x,y\in X$ with $d(x,y)=\sigma$, we choose a ball $B$ centered in the center of a geodesic connecting $x$ and $y$. Then we can connect $x$ and the closest of the points $x'$ and $x''$ (say, $x'$) with a Harnack chain away from $x''$ because $d(x',x'')>\sigma/9$, i.e., we place a sequence of overlapping equally sized balls along a geodesic connecting these points and apply the Harnack inequality a controlled number of times. This gives a controlled bound on $G(x,x'')$ and using a similar Harnack chain connecting $x''$ and $y$ (which now only has to avoid $x$) we get the upper bound on $G(x,y)$.
The ball volumes can be compared using the doubling property.

\end{proof}
\end{proposition}
Note that we could now easily get more explicit growth estimates for $G$ near the pole by comparison with the explicit estimates for local Green functions, but for the following these coarse estimates are sufficient.

\subsubsection{Relative Maximum Principles}
We start with an elementary comparison result which already visualizes the primary effect of coercivity: with the same boundary conditions, $\L$-harmonic functions \emph{sag} compared to $\L_{-\epsilon}$-harmonic functions. Imagine a rope\footnote{This analogy is quite precise, a suspended rope forms a hyperbolic cosine, a solution of the shifted one-dimensional Laplace equation.} or a rubber blanket where you apply less and less tension. To see this, we start with a sufficiently curved function that fits in between.

\begin{lemma}
There is a constant $m(\sigma, N, C_P, C_D, k, \epsilon)>0$ such that on any ball $B_{\sigma}(x_0)$, we can find a continuous function $f\in\Floc(B_{\sigma}(x_0))$ with
\[ f(x_0)\le-m\,,\quad f|_{\del B_{\sigma}(x_0)}\ge 0\quad\text{and}\quad \L f\ge-1\text{ on $B_\sigma(x_0)$.}\]
\begin{proof}
We consider the Dirichlet Green function $G_{B}$ for $\L$ on $B:=B_{\sigma}(x_0)$, then $G_{B}(\chi_{B})(x_0)\ge G_{B}(\chi_{B_{\sigma}(x_0)\setminus B_{\sigma/2}(x_0)})(x_0)\ge m$ for some universal constant $m>0$ by doubling and \thmref{prop:bgreen}.
Hence, the function $f:=-G_{B}(\chi_{B})$ satisfies $\L f\equiv-1$ on $B$, $f|_{\del B}\equiv 0$ and $f(x_0)\le -m$.

\end{proof}
\end{lemma}

\begin{proposition}[relative maximum principle, local version]
\label{qmploc}
Assume that we have two functions $u$ and $\bar u>0$ on a ball $B_\sigma=B_\sigma(x_0)$, $u$ $\L$-subharmonic ($\L u\le0$) and $\bar u$ $\L_{-\epsilon}$-harmonic on $B_\sigma$ with $\bar u|_{\del B_{\sigma/4}}\ge u|_{\del B_{\sigma/4}}$.
Then there is a constant  $\tilde\eta=\tilde\eta(\sigma, N, C_P, C_D, k, \epsilon) \in (0,1)$ such that
\[ u(x_0)\le \tilde\eta\,\bar u(x_0)\,.\]
\begin{proof}
Consider the function $h(z):=\bar u(z)+\epsilon f(z)\inf_{w\in B_{\sigma/4}}\bar u(w)-u(z)$ on $B_{\sigma/4}$ where $f$ is the function from the preceding lemma.
On $B_{\sigma/4}$, we have $\L\bar u=\L_{-\epsilon}\bar u+\epsilon\bar u\ge \epsilon\bar u$, $\L f>-1$ and therefore
\[\L h(z)>\epsilon\left(\bar u(z)-\inf_{w\in B_{\sigma/4}}\bar u(w)\right) \ge 0\,.\]
The boundary condition together with $f|_{\del B_{\sigma/4}}\ge 0$ yields $h|_{\del B_{\sigma/4}}\ge 0$ such that by the \thmnameref{prop:min0}
$h\ge 0$ in the interior of $B_{\sigma/4}$ and especially $h(x_0)\ge0$; with $f(x_0)\le-m$ we have
\[ u(x_0)\le\bar u(x_0)-m\epsilon\inf_{w\in B_{\sigma/4}}\bar u(w)\le(1-m \epsilon H^{-1})\bar u(x_0)\,. \]
The harmonicity on the full ball $B_\sigma$ was used only in the last step to apply the Harnack inequality.
\end{proof}
\end{proposition}

Applied globally, this describes the relative growth of  $\L$-harmonic versus $\L_{-\epsilon}$-harmonic functions.

\begin{proposition}[relative maximum principle, global version]
\label{qmpglob}
There is a constant $\eta=\eta(\sigma, N, C_P, C_D, k, \epsilon)\in(0,1)$ such that the following holds:

Assume we have two functions $u$ and $\bar u>0$ defined on $B_{r+\sigma}(x)$ for some $x\in X$ and $r\ge \sigma$, $u$ $\L$-subharmonic and $\bar u$ $\L_{-\epsilon}$-harmonic on $B_{r+\sigma}(x)$ and $\bar u|_{\del B_r(x)}\ge u|_{\del B_r(x)}$. Then
\[ u(x)\le \eta^r \bar u(x)\,.\]
\begin{proof}
For integer multiples $r$ of $\sigma/4$, this is proven by inductively applying \hyperref[qmploc]{the local version} along a chain of intersecting balls of length proportional to $\sigma/4$. On each of them we may apply the Harnack inequality as described in the proof of \autoref{prop:bgreen} above, so that we can choose $\eta$ as a function of $\sigma$, $H$ and $\tilde\eta$.
\end{proof}
\end{proposition}

For Green functions, we get the following variants. We do not use them in the following arguments, but they are worth being mentioned since they give some non-trivial constraints on the Green functions without using Gromov hyperbolicity.

\begin{corollary}[exponentially stronger decay]
\label{prop:expstrdec}
There are universal constants $A(\sigma, N, C_P, C_D, k, \epsilon)>0$ and $\alpha_1(\sigma, N, C_P, C_D, k, \epsilon)>0$ such that
\[ G(x,y)\le A\e^{-\alpha_1 d(x,y)}G_{-\epsilon}(x,y)\quad\forall x,y\in X\,. \]
\end{corollary}
From this we get the following growth estimate for Green functions using the resolvent equation, cf.~\cite[Corollaire~4.6]{Anc90}:
\begin{proposition}[exponential decay]%TODO called 4.7 in Scal Splittings I arXiv:2012.12223 [math.DG]
\label{prop:expdec}
For suitable constants $B(\sigma, N, C_P, C_D, k, \epsilon)>0$ and $\alpha_2(\sigma, N, C_P, C_D, k, \epsilon)>0$ we have
\[ G(x,y)\le B \e^{-\alpha_2 d(x,y)}\quad\text{for }d(x,y)> 2\sigma\,.\]
%In particular, the Green function of a symmetric operator decays exponentially.
\begin{proof}
Let $x'\in X$ such that $d(x,x')=\sigma$. Then employing the symmetry of the Green function and the \thmnameref{resolvent} we have
\begin{multline*}
	G(x,y)^2\le G(x,y)\,H^3\,G(y,x')
\le \frac{H^5}{\vol(B_{\sigma/2}(y))}\int_{B_{\sigma/2}(y)} G(x,z)\,G(z,x')\,\d\mu(z)\\
    \le c\int_X G(x,z)\,G(z,x')\,\d\mu(z)
	\overset{\eqref{eq:resineq1},\eqref{eq:resineq2}}{\le} \frac{c}{\epsilon}G_{-\epsilon}(x,x')\le\frac{cc_1^\epsilon}{\epsilon}
\end{multline*}
where we applied the \thmnameref{prop:bgreen} to $\L_{-\epsilon}$ which itself satisfies the assumptions, but with a weaker constant $c_1^\epsilon$.

For the very same reasons, we can do all of the above with $G_{-\epsilon}$ instead of $G$ to get the uniform boundedness of $G_{-\epsilon}(x,y)^2$, again with slightly worse constants. Combined with \thmref{prop:expstrdec} we have proved the assertion.
\end{proof}
\end{proposition}

This result does not use Gromov hyperbolicity and holds also e.g.\ in Euclidean space. One may remember that the familiar Laplacian's Green function on $\R^n$ only decays with $|x-y|^{-(n-2)}$, but this is not in conflict to the result above because the Euclidean Laplacian is not coercive.

\section{Potential Theory}\label{sec:pottheo}
Before proving the main results, we will introduce central potential theoretic and geometric concepts in this and the following section, in particular two adapted notions of boundary at infinity classifying globally defined harmonic functions (Martin boundary) and geodesic rays (Gromov boundary). In \autoref{sec:ident}, they will turn out to be homeomorphic.
\subsection{Axioms and Potentials}\label{sec:axpot}
We will employ several constructions that are best described in the language of abstract potential theory. Here we will briefly describe central concepts. More background information and details can be found in standard literature on the subject such as \cite{Bau66, Bre67, Hel69, CC72, BH86}.

A common set of basic axioms for harmonic functions on a locally compact space $X$ is the following:
\begin{enumerate}
  \item The harmonic functions associated with any open subset of $X$ constitute a sheaf of continuous functions. \textbf{(Sheaf Property)}
  \item The topology of $X$ has a basis consisting of relatively compact open sets $U$ where each continuous function $f$ on the boundary $\del U$ has a unique continuous extension to $\bar U$ that is harmonic on $U$ and positive if $f$ is positive. \textbf{(Basis of Regular Sets)}
  \item The (pointwise) limit of an increasing sequence of harmonic functions on a connected open set $U\subset X$ is either $\equiv+\infty$ or harmonic. \textbf{(Brelot Convergence Axiom)}
\end{enumerate}

These properties are in fact satisfied by the $\L$-harmonic functions with respect to a coercive Schrödinger operator $\L$ in our setting:
\begin{enumerate}
\item The Sheaf Property follows easily from our definition of harmonic functions, as soon as there are partitions of unity consisting of continuous functions in $\F$. Their existence in the case of regular Dirichlet forms is proven in \cite[Theorem~3.1]{AMR97}. %product in $\F$: \cite[Corollary 4.15]{MR92}
\item A basis of regular sets for $\L^0$ is given in \autoref{sec:strreglocbounds}. In \autoref{sec:locbounds}, we have seen that these sets of a sufficiently small diameter are still regular for $\L$, giving a basis in this case as well. Positivity is ensured by the \thmnameref{prop:min0}.
\item Given the first two axioms, Brelot's convergence axiom is equivalent to the Harnack inequality we proved in \autoref{sec:locbounds} \cite{LW65}.
\end{enumerate}

In this axiomatic setting, \textbf{$\bm{\L}$-superharmonic} functions can be characterized as lower semi-continuous functions $u$ with values in $(-\infty,\infty]$ such that $u$ is finite on a dense set and $u\ge v_U$ for every regular set $U$, where $v_U$ is the solution of $\L v_U=0$ on $U$ with boundary condition $v_U|_{\del U}\equiv u|_{\del U}$. This is a generalization of the notion of superharmonic functions introduced in \autoref{sec:dirichlet},
and slightly abusing notation, we will write $\L u\ge0$ for a superharmonic function $u$.

If they exist, as in the case of \emph{coercive} Schrödinger operators, Green functions are superharmonic, and they usually satisfy a (generalized) vanishing condition at the boundary/infinity. This additional property is captured in the notion of an \emph{\textbf{$\bm{\L}$-potential}}, which is an $\L$-superharmonic function $p>0$ such that there is \emph{no} positive $\L$-harmonic function $h\in\F_\loc(X)$ with $p \ge h >0$. Green functions are not only examples for potentials, but also their basic building blocks:

\begin{theorem}[integral representation of potentials]\cite[22.]{Her62}, \cite[Thm.\,6.18]{Hel69}
If $\L$ has a Green function $G$, every $\L$-potential $p$ on $X$ can be represented by a unique (positive) Radon measure $\mu_p$ as
\[p(x)=G(\mu_p)(x):=\int_X G(x,y)\,\d\mu_p(y).\]
The \emph{support} of $p$ (i.e., the complement of the largest open set where $p$ is $\L$-harmonic) equals the support of $\mu_p$.
\end{theorem}

The fact that every positive $\L$-superharmonic function is uniquely representable as the sum of an $\L$-potential and a positive $\L$-harmonic function is known as the \emph{Riesz decomposition}. Hence, to get an integral representation for positive \emph{$\L$-superharmonic} functions, only the $\L$-harmonic part is left. This can be taken care of in an abstract way as described in the following.

\subsection{Martin Boundary}
\label{sec:martin}

In the same spirit as $\L$-potentials can be represented by Green functions integrated over Radon measures on $X$, there is a representation of positive $\L$-harmonic functions by Radon measures on an abstractly associated space $\delM(X,\L)$, the \emph{Martin boundary}. A priori, this is a purely potential theoretic object and there might be additional peculiarities in the representation, but as a central result of \autoref{sec:ident} we will see that the Martin boundary for a coercive Schrödinger form on a Gromov hyperbolic space can be canonically identified with the Gromov boundary.

Here, we sketch the basic notions from Martin theory, see~\cite[§11.4]{CC72}, \cite{Anc90}, \cite[7.1]{Pin95} or \cite[Ch.\,12]{Hel69} for detailed proofs.

\begin{definition}[Martin Boundary] In a non-compact space $X$ satisfying the axioms from the preceding subsection where a Green function $G: X\times X \ra (0,\infty]$ exists, let $o \in X$ be a basepoint. We consider the set $S$ of sequences $s$ of points $x_i \in X$, $i=1,2,\dots$, such that
\begin{itemize}
    \item $s$ has no accumulation point in $X$, and
    \item $K_{x_i}:= \frac{G(\cdot,x_i)}{G(o,x_i)} \overset{i\to\infty}{\longrightarrow} K_s$ locally uniformly, for some function $K_s:X\to(0,\infty)$.
\end{itemize}
As a set, the \textbf{Martin boundary} $\delM(X,\L)$  is the quotient of $S$ modulo the relation $s \sim s'$ if $K_s \equiv K_{s'}$. For any $\zeta\in \delM(X,\L)$, this function is written $K_\zeta$ and is called a \textbf{Martin function}.
\end{definition}

The Martin boundary does not depend on the choice of the basepoint $o$.  The Harnack inequality and elliptic theory show that each $K_\zeta \in \delM  (X,\L)$ is a positive $\L$-harmonic function on $X$. This also shows
that the convex set $S_\L(X)$ of positive $\L$-harmonic functions $u$ on $X$ with $u(o)=1$ is compact in the topology of local uniform convergence. In turn, $\delM  (X,\L)$ is a compact
subset of $S_\L(X)$.

We topologize the space $\overline{X}^\mathrm{M}:= X \cup \delM(X,\L)$ using the topology of local uniform convergence on the space of associated Martin functions $\{K_y\:|\:y\in\overline{X}^\mathrm{M}\}$, which can be canonically identified with $\overline{X}^\mathrm{M}$.
Then the usual topology is induced on $X\subset\overline{X}^\mathrm{M}$, $\delM(X,\L)$ is closed and $\overline{X}^\mathrm{M}$ is compact. The space $\overline{X}^\mathrm{M}$ is called the \textbf{Martin compactification}. It is easy to see that $\overline{X}^\mathrm{M}$ is metrizable, see \cite[Ch.\,I.7]{BJ06} or \cite[Ch.\,12]{Hel69} for further details.

To motivate the idea of Martin integrals, we recall the following classical result, see e.g.\ \cite[Ch.\,6]{Cho69}:

\begin{proposition}[Minkowski's theorem]\label{mit} Each point in a convex set  $K \subset \R^n$ is a
convex combination of the extremal points of $K$.
\end{proposition}

\begin{figure}[h]
\centering
\includegraphics[width=0.84\textwidth]{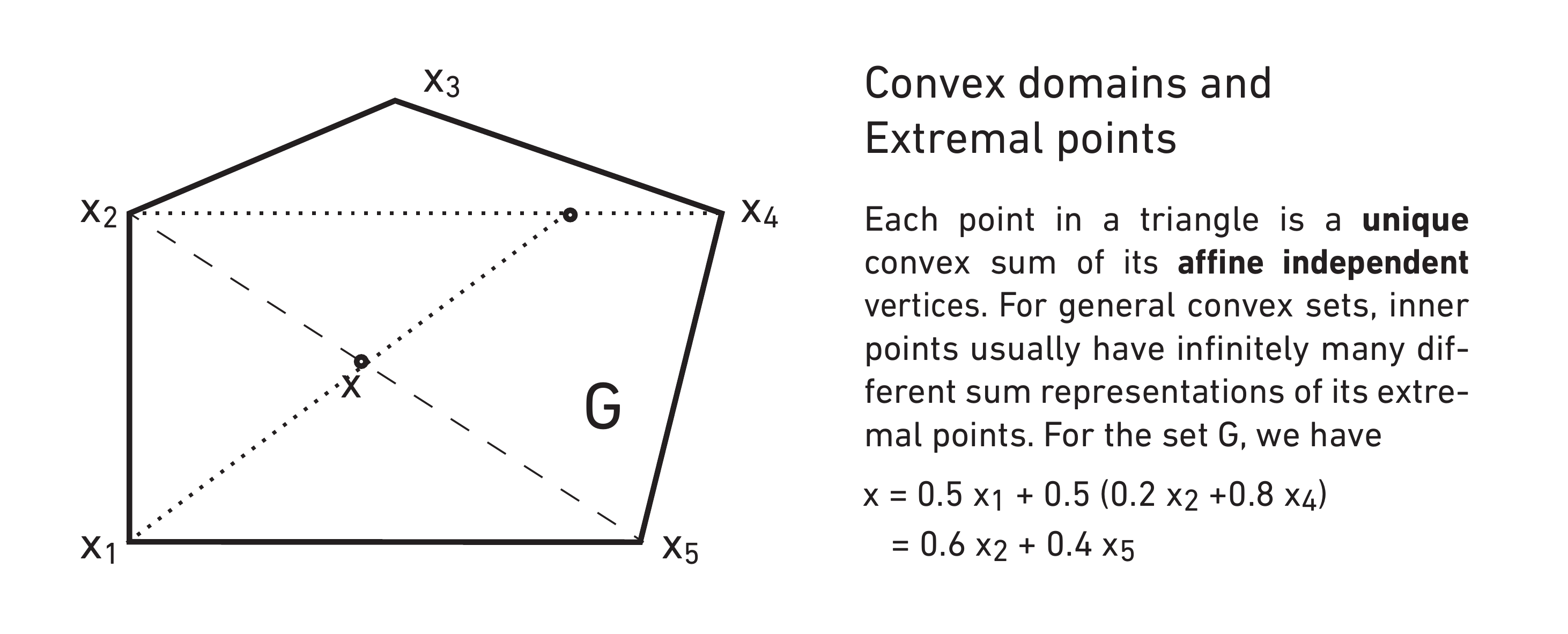}
\caption{Minkowski's Theorem in $\R^2$.}
\end{figure}

The \textbf{Martin integral representation} is essentially an extension of this result to the case of the convex set $S_\L(X)$. The extremal elements of $S_\L(X)$ form a subset $\delM^0 (X,\L) \subset \delM(X,\L)$ of the Martin boundary which can be seen as the vertices of an infinite dimensional simplex spanning $S_\L(X)$. A positive solution $u$ of $\L u = 0$ on $X$ with $u(o)=1$ is extremal if and only if $u$ is a \emph{minimal} solution in the following sense: for any other solution $v>0$,  $v \le u$, there is a constant $c>0$ such that $v \equiv c \cdot u$. Therefore, $\delM^0 (X,\L) \subset \delM (X,\L)$ is also called the \textbf{minimal Martin boundary}.

The Choquet integral representations in \cite[Ch.\,6]{Cho69} give the following general version of the Martin representation theorem, see \cite[7.1]{Pin95}:

\begin{proposition}[Martin integral representation]\label{thm:mrt} For any positive solution $u$ of $\L u = 0$ on $X$, there is a
unique Radon measure $\mu_{u}$ on $\delM^0 (X,\L)$, sometimes called \emph{Martin measure of $u$}, so that
\[u(x)  =\int_{\delM^0 (X,\L)} K_\zeta(x) \, \d \mu_u(\zeta)\,.\]
Conversely, for any Radon measure $\mu$ on $\delM^0 (X,\L)$,
\[u_\mu(x)  =\int_{\delM^0 (X,\L)} K_\zeta(x) \, \d \mu(\zeta)\]
defines a positive solution of $\L u_\mu = 0$ on $X$.
\end{proposition}

Although this already looks like a classical contour integral, with $K$ generalizing the Poisson kernel, the result is not yet truly satisfactory. Unlike the classical case, the boundary $\delM^0 (X,\L)$ depends not only on the underlying space, but also on the analysis of the operator $\L$. A natural question is whether one could get rid of this dependence. In general the answer is no, as is shown e.g.\ in \cite{Anc12}, where Ancona describes non-uniform Euclidean cones with only one topological point at infinity but with uncountably many minimal Martin boundary points at infinity.
However, we will see that in our case of coercive Schrödinger operators on Gromov hyperbolic spaces, the Martin boundary can be completely described in terms of the geometric Gromov boundary. This is a remarkable particularity not even valid for such simple spaces as in \thmref{ex:pro}.

\subsection{Balayage}
\label{sec:balayage}
Later we want to control $\L$-superharmonic functions along $\Phi$-chains. Here we shift the part of the defining measure supported in $X \setminus  U$ onto $\del U$ without changing the function on $U$, for an unbounded open set $U$. This strategy is called sweeping or, due to its French origin (Poincar\'{e}, Cartan), \emph{\textbf{balayage}}.

Concretely, for an $\L$-superharmonic function $u\ge0$  on $X$ and a subset $A\subset X$ we define
\[\red_u^A:= \inf \{ v \ge0\:|\:v \mbox{ is $\L$-superharmonic on } X \mbox{ with } v \ge u \mbox{ on } A \}\,.\]
This is called the \emph{\textbf{reduit}} (reduced).
It enjoys the following properties which we will need later:
\begin{itemize}
\item $\red_u^A$ is $\L$-harmonic outside $\bar A$ and equal to $u$ on $A$.
\item The reduit is always $\L$-superharmonic.
\item If $A$ is relatively compact, $\red_u^A$ is an $\L$-potential.
\item $\red_{\lambda u}^{A}=\lambda \red_{u}^{A}$ for a constant $\lambda\ge0$.
\item $\red_{u+v}^{A}=\red_{u}^{A}+\red_{v}^{A}$ for functions $u$, $v$.
\item $\red_u^{A\cup B}\le \red_u^A+\red_u^B$ for sets $A,B\subset X$.
\item Since $\L$ is self-adjoint and hence $G$ symmetric, we have $\red_{G(\cdot,y)}^A(x)=\red_{G(x,\cdot)}^A(y)$.
\item If $G(\mu)$ is an $\L$-potential, $\red_{G(\mu)}^A(x)=\int_X\red_{G(\cdot,y)}^A(x)\,\d\mu(y)$ for any $x\notin\bar{A}$ \cite[Théorème 22.4]{Her62}.
\end{itemize}

For general sets $A$, it may happen that $\red_u^A$ is not lower semi-continuous, but it always admits a canonical regularization $\widehat{\red}_u^A$, the \emph{balayée} (swept), defined as the maximal lower semi-continuous function $\le\red_u^A$. For open sets $A$ or in general outside $\bar{A}$ the two concepts coincide.
We can even recover the classical Perron solution $u$ of the Dirichlet problem on a ball $B$ with continuous positive boundary value $f$ as
\[ u(x)=\red_f^{\del B}(x)\,.\]
Note that this allows us to solve Dirichlet problems on arbitrarily large regular bounded sets as soon as there is a Green function, e.g.\ in the case of coercive Schrödinger forms.

\subsection{A Global Maximum Principle}
As an application of the axiomatic approach, we note two versatile variants of the maximum principle for later use.
\begin{theorem}[global maximum principle]\cite[p.\,429]{Her62}\label{gmp}
\begin{enumerate}
\item If $u$ is $\L$-superharmonic on an open set $V\subset X$, $u\ge0$ on $\del V$, and there is an $\L$-potential $p$ such that $u\ge-p$, then $u\ge0$ on $V$.
\item Let $p$ an $\L$-potential, $\L$-harmonic on an open set $V$ and locally upper bounded near every point of $\del V$. If $u\ge p$ on $\del V$ for some positive $\L$-superharmonic function $u$, then $u\ge p$ in all of $V$.
\end{enumerate}
\begin{proof}
For (i), note that the function $\bar u$ defined as $\min(u,0)$ on $V$ and $0$ on $X\setminus V$ is $\L$-superharmonic and $\ge-p$. Now the supremum of the family $\{\text{$\L$-subharmonic functions $\le \bar u$}\}$ is $\L$-harmonic,
$\ge -p$, and $\le0$, hence by the definition of $\L$-potentials it is 0 which implies $u\ge0$.
\medskip

(ii) follows from (i) by considering the function $u-p$.
\end{proof}
\end{theorem}

\section{Gromov Hyperbolic Spaces}\label{sec:gromhyp}
This additional geometric condition on our geodesic metric space $X$ is complementary to the bounded geometry assumptions in the sense that it has no local impact, but rather limits all efficient communication between two points to a region near the shortest path between them.

\begin{figure}[h]
  \centering
  \includegraphics[width=7cm]{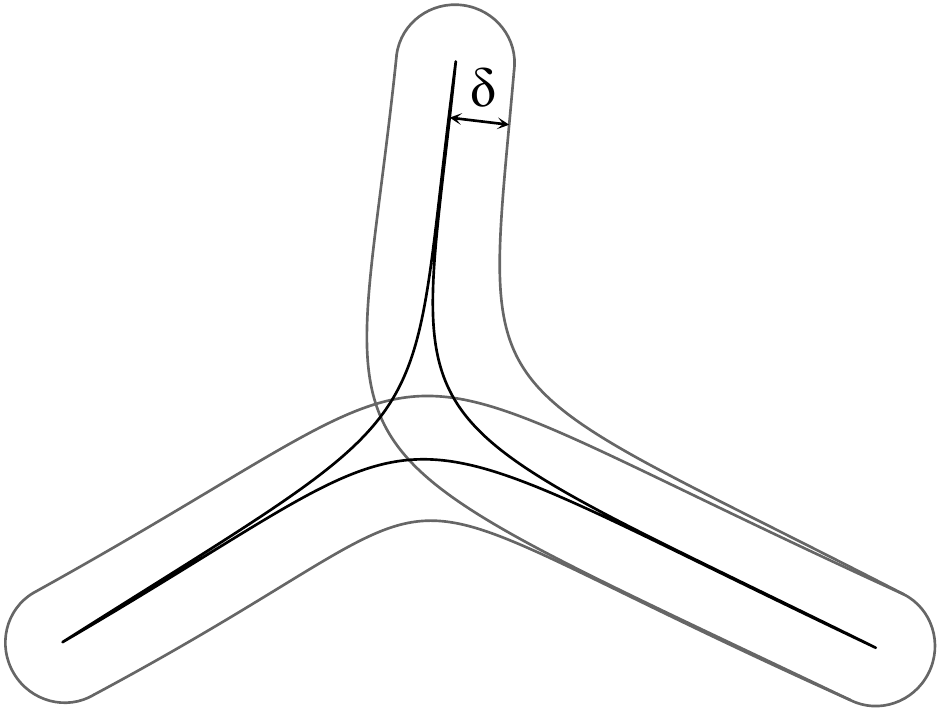}
  \caption{A geodesic triangle in a $\delta$-hyperbolic space.}\label{fig:deltathin}
\end{figure}

\begin{definition}[Gromov hyperbolicity]\label{def:hyp}
A geodesic metric space is \emph{Gromov hyperbolic} or, quantitatively, \textbf{$\bm{\delta}$-hy\-per\-bo\-lic}, if there is a $\delta>0$ such that each point on the edge of any geodesic triangle\footnote{Here, \enquote{geodesic} is meant in the metric space sense as \enquote{length-minimizing geodesic}.} is within $\delta$-distance of one of the other two edges. See \autoref{fig:deltathin}.
\end{definition}

$\delta$-hyperbolicity was introduced for Cayley graphs of finitely generated groups by Gromov in \cite{Gro87}. Basic examples are:

\begin{examples}\leavevmode
\begin{enumerate}
  \item Bounded spaces are trivially $\delta$-hyperbolic for $\delta$ equal to the diameter.
  \item Hyperbolic space $\mathbb{H}^n$ is $\delta$-hyperbolic for $\delta=\ln 3$ \cite[2.23]{Gd90}. More generally, $\mathrm{CAT}(-1)$-spaces are $\ln3$-hyperbolic, as directly seen from the definition \cite[Proposition~1.2]{BH99}.
  \item \emph{Cartan--Hadamard manifolds}, i.e., complete simply-connected manifolds with sectional curvature bounded from above by a negative constant, are $\mathrm{CAT}(\kappa)$-spaces for some $\kappa<0$ and hence Gromov hyperbolic \cite[§3.2]{Gd90}. If the sectional curvature is additionally bounded from below, they have also bounded geometry.
  \item Particularly simple examples that can already illustrate many phenomena are graphs. By a (metric) \textbf{\emph{graph}} we mean the geodesic metric space obtained by gluing intervals of a certain length between elements of a potentially infinite set of \emph{vertices}, such that each vertex is adjacent to a finite number of these \emph{edges}. Bounded geometry (and finite-dimensionality) of manifolds is analogous to a global upper bound on the number of adjacent edges for each vertex and a lower bound on their length.

      \emph{Trees} (simply connected graphs) are precisely the 0-hyperbolic graphs, hence for $\delta$-hyperbolic graphs, $\delta$ can be seen as a measure of deviation from being a tree.
  \item A finitely generated group is called hyperbolic, if the Cayley graph with respect to some (and hence any, see \cite[Examples~I.8.17 (2) and (3)]{BH99}) finite set of generators with edge length 1 is Gromov hyperbolic. Here, the easiest example are free groups (corresponding to trees, for canonical generators).
  \item The universal covering of a closed manifold $X$ of negative sectional curvature satisfies the requirements in (iii). By the Švarc-Milnor Lemma \cite[Prop.~I.8.19]{BH99}, any Cayley graph of the fundamental group $\pi_1(X)$ is quasi-isometric to $X$ and hence also Gromov hyperbolic.

      Here, a map $f:X\to Y$ between metric spaces $(X,d_X)$, $(Y,d_Y)$ is a $(\lambda, S)$-\textbf{\emph{quasi-isometry}} if there are constants $\lambda\ge1$, $S\ge0$ such that
      \[
        \lambda^{-1}d_X(x,x')-S\le d_Y(f(x),f(x')) \le \lambda d_X(x,x') + S\quad\forall x,x'\in X
      \]
      and every point $y\in Y$ has distance at most $S$ to a point in the image of $f$. Without the last condition, $f$ is a $(\lambda, S)$-\textbf{\emph{quasi-isometric embedding}}. If such a quasi-isometric map exists, $X$ and $Y$ are \emph{quasi-isometric}.

      For invariance of Gromov hyperbolicity under quasi-isometries see e.g.\ \cite[1.3.1]{BS07}.
  \item A non-example: Gromov hyperbolicity can be more demanding than constant negative sectional curvature alone. Consider a $\Z^2$-covering of a Riemann surface of genus $\ge2$ equipped with a metric of constant negative sectional curvature, as in \autoref{fig:hypnothyp}.
      %\begin{wrapfigure}{R}{0.5\textwidth}
      \begin{figure}
        \centering
        \includegraphics[width=0.7\textwidth]{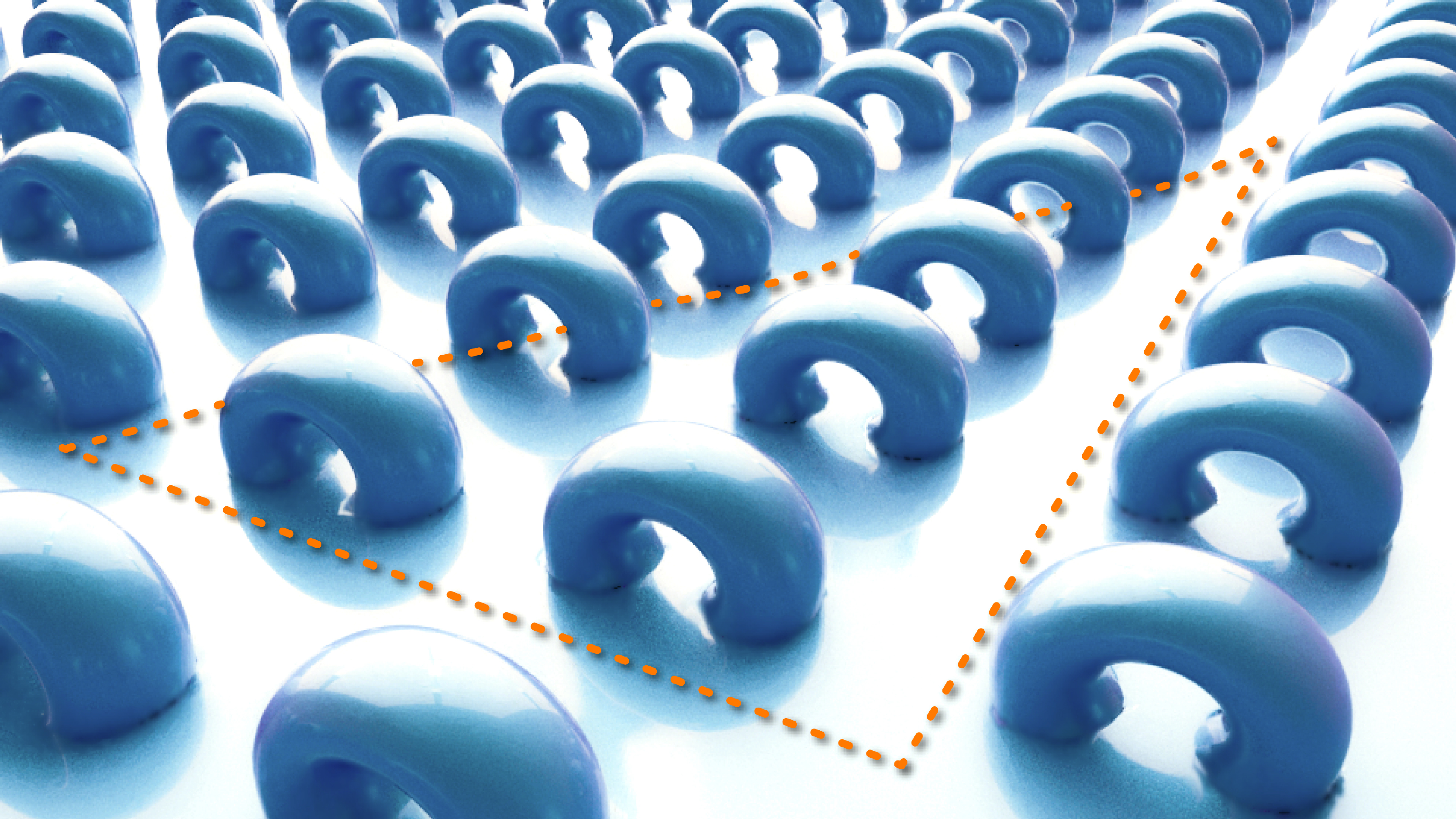}
        \caption{We can get complete 2-manifolds with \emph{constant negative sectional curvature} by periodically attaching handles to the flat plane. Nevertheless, this space is \emph{not} $\delta$-hyperbolic since it coarsely looks like the plane and the typical geodesic triangles are nearly Euclidean.}\label{fig:hypnothyp}
      \end{figure}
      %\end{wrapfigure}
      This is quasi-isometric to Euclidean $\R^2$ and hence not Gromov hyperbolic.
\end{enumerate}
\end{examples}

The idea alluded to above that there is no efficient path between two points that is far away from the shortest path can be formalized as follows:
\begin{proposition}[stability of geodesics]\cite[Theorem~1.3.2]{BS07}, \cite{Bon96}\label{prop:geostab}
In a $\delta$-hyperbolic geodesic metric space $X$, for every $\lambda\ge1$ and $S>0$ there is an $H=H(\delta, \lambda, S)>0$ such that each two $(\lambda, S)$-\emph{quasi-geodesics}\footnote{$(\lambda,S)$-quasi-isometric embeddings of a compact interval} $x\leadsto y$ with same start- and endpoint $x,y\in X$ have Hausdorff distance at most $H$. In fact, this \emph{geodesic stability} is equivalent to Gromov hyperbolicity.
\end{proposition}

There are several definitions of Gromov hyperbolicity highlighting different aspects, see e.g.\ \cite[Chapter~1]{BS07}, \cite[Chapter~III.H]{BH99} or \cite[Chapitre~2]{Gd90}. We will only mention one other version using the \textbf{\emph{Gromov product}} on $X$ defined as
\[
    (x|y)_z := \frac12\left(d(z,x)+d(z,y)-d(x,y)\right)\quad\text{for $x,y,z\in X$.}
\]
This version works on arbitrary metric spaces.
\begin{definition}[hyperbolicity via Gromov product]
A metric space $X$ is \emph{Gromov hyperbolic} if there is a $\delta'\ge 0$ such that for any four points $x,y,z,w\in X$,
\begin{equation}\label{eq:hypdef}
    (x|y)_z\ge\min\{(x|w)_z,(w|y)_z\}-\delta'\,.
\end{equation}
This implies $4\delta'$-hyperbolicity and is implied by $\delta'/8$-hyperbolicity as defined above, see \cite[Proposition~2.21]{Gd90}. For convenience (and since the precise constant never really matters) we assume from now on every $\delta$-hyperbolic space to satisfy \eqref{eq:hypdef} with $\delta'=\delta$.
\end{definition}

\begin{figure}
  \centering
  \includegraphics[width=7cm]{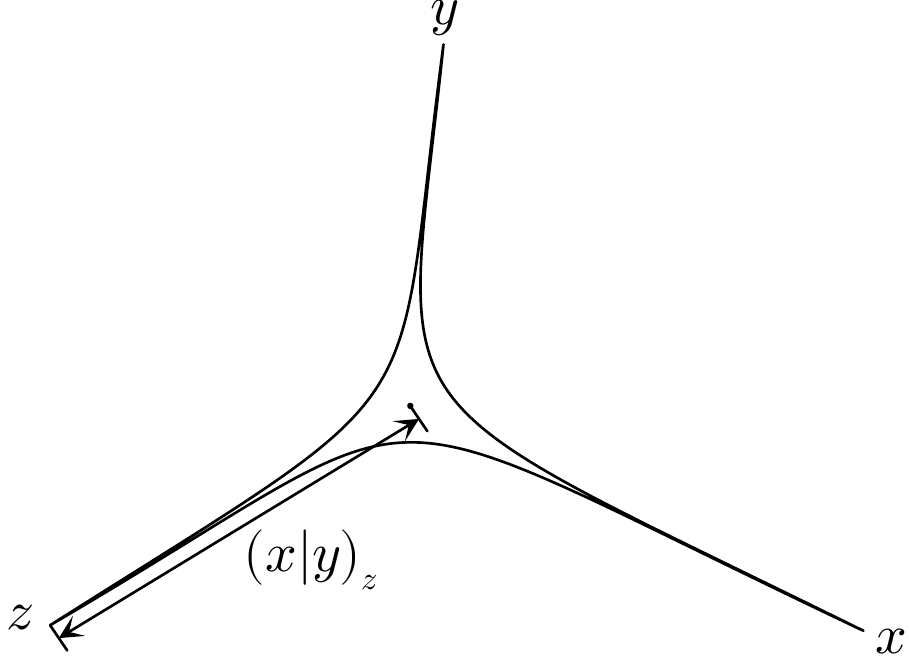}
  \caption{Approximate value of the Gromov product in a $\delta$-hyperbolic space.}\label{fig:gromovproduct}
\end{figure}

An intuitive interpretation of the Gromov product in Gromov hyperbolic spaces is given by the following estimate (see also \autoref{fig:gromovproduct}):
\begin{lemma}[Gromov product as distance to a geodesic]\label{lem:productasdistance}\cite[Lemme~2.17]{Gd90}
In a $\delta$-hyperbolic geodesic metric space $X$, let $\gamma:x\leadsto y$ be a geodesic and $z\in X$. Then
\[
    (x|y)_z \le \dist(z,\gamma)\le (x|y)_z + 4\delta\,.
\]
\end{lemma}

\subsection{Gromov Boundary}\label{sec:gromovboundary}
In Gromov hyperbolic spaces, the structure at infinity can be encoded in an ideal boundary. We outline the construction presented in \cite[Section~2.2]{BS07}, where more details can be found.
\begin{definition}[Gromov boundary]
In a Gromov hyperbolic metric space $X$ with basepoint $o\in X$, a sequence $(x_i)$ in $x$ \emph{converges at infinity} if $(x_i|x_j)_o\overset{i,j\to\infty}{\longrightarrow}\infty$. We say two sequences $(x_i)$, $(y_i)$ convergent at infinity are \emph{equivalent} if $(x_i|y_i)_o\overset{i\to\infty}{\longrightarrow}\infty$. The \textbf{\emph{Gromov boundary}} $\delG X$ of $X$ is defined as the set of equivalence classes of sequences convergent at infinity.
\end{definition}

Note that $\delG X$ is independent of the chosen basepoint because $|(x|y)_o-(x|y)_{o'}|\le d(o,o')$ for a different basepoint $o'\in X$.

To define a topology on $\overline{X}^{\mathrm{G}} := X \dot\cup \delG X$, we first extend the Gromov product to the boundary. We formally identify points in $X$ with constant sequences. Then for $a,b\in\overline{X}^{\mathrm{G}}$ and $z\in X$, we can set
\[
    (a|b)_z:=\inf_{\substack{(x_i)\in a\\(y_i)\in b}}\liminf_{i\to\infty}(x_i|y_i)_z\,.
\]
A topology on $\overline{X}^{\mathrm{G}}$ is defined by considering $X\hookrightarrow \overline{X}^{\mathrm{G}}$ as an embedding, choosing a basepoint $o\in X$ and declaring the sets
\[
    \W_\rho^o(\xi):=\{a\in\overline{X}^{\mathrm{G}}\:|\:(a|\xi)_o>\rho\}\quad\text{for $\rho\ge0$}
\]
to be a neighborhood basis for $\xi\in\delG X$. As above, this topology is independent of $o$.

The sets $\W_\rho^o(\xi)\cap\delG X$ are in fact the open balls of radius $\e^{-\rho}$ with respect to
\[
    d_o(\xi,\eta):=\e^{-(\xi|\eta)_o}\,,
\]
which is a \emph{quasi-metric} on $\delG X$, i.e., an ultrametric triangle inequality holds only up to a constant $Q=Q(\delta)=\e^{\delta}$,
\[
    d_o(\xi,\zeta)\le Q\max\{d_o(\xi,\eta),d_o(\eta,\zeta)\} \text{ for $\xi,\eta,\zeta\in\delG X$,}
\]
while all other properties of metrics are still satisfied. This quasi-metric triangle inequality follows from an extension of \eqref{eq:hypdef} to the boundary \cite[Lemma~2.2.2\,(2)]{BS07}, while the other properties are obvious.

\begin{remarks}\leavevmode\label{rem:hypapprox}
\begin{enumerate}
  \item If the Gromov hyperbolic space $X$ is proper (balls are relatively compact) and geodesic (e.g., a complete manifold), it is sufficient to consider sequences on geodesic rays emanating from the basepoint $o$. One could also consider equivalence classes consisting of geodesic rays with finite Hausdorff distance \cite[2.4.2]{BS07}. Hence, any geodesic ray has a well-defined endpoint in $\delG X$, and we may use the notation $\gamma:x\leadsto\xi$ for a geodesic ray $\gamma$ starting in $x\in X$ with endpoint $\xi\in\delG X$.
  \item For a proper geodesic Gromov hyperbolic space $X$, $\overline{X}^{\mathrm{G}}$ is in fact a compactification, i.e., $X\subset\overline{X}^{\mathrm{G}}$ is open and dense and $\overline{X}^{\mathrm{G}}$ is compact \cite[Proposition~III.H.3.7]{BH99}.
  \item For quasi-metrics with constant $Q\le2$, there is a canonical procedure to construct a bi-Lipschitz equivalent metric \cite[2.2.2]{BS07}. This can be applied to the quasi-metrics $d_o^\epsilon$ for sufficiently small $\epsilon$ to get a family of metrics on the boundary (called \emph{visual metrics}), but we prefer to work with the canonical quasi-metric instead.
  \item There are \emph{many} hyperbolic spaces in the following sense: for each bounded metric space $Z$, there is a \emph{hyperbolic approximation} of $Z$: this is a Gromov hyperbolic graph $X$ with basepoint $o$ such that the boundary $\delG X$ equipped with the quasi-metric $d_o$ is bi-Lipschitz equivalent to $Z$ \cite[Theorem~6.4.1]{BS07}. $X$ is proper if and only if $Z$ is compact \cite[6.4.3]{BS07}.
  \item The graph in the preceding remark can be approximated by a manifold of any dimension $n\ge2$. To this end, represent each vertex by an $n$-sphere and for each edge connecting two vertices, form a connected sum of the representing spheres. If there is a global upper bound on the number of edges adjacent to any vertex, the metric on this manifold can be arranged to have bounded geometry.
\end{enumerate}
\end{remarks}

\subsection{\texorpdfstring{\textit{Φ}}{Phi}-Chains}

\paragraph{Gromov Hyperbolicity} To exploit Gromov hyperbolicity analytically, Ancona introduced the concept of \emph{$\Phi$-chains}. Just as the large-scale notion of hyperbolicity complements small-scale bounded geometry, $\Phi$-chains act complementary to Harnack chains. This is best seen in the discussion preceding \thmref{prop:expg} where Harnack chains give basic estimates which, however, weaken the overall control, while $\Phi$-chains can be used to recover the apparently lost details.

\begin{figure}
  \centering
  \includegraphics[width=.76\textwidth]{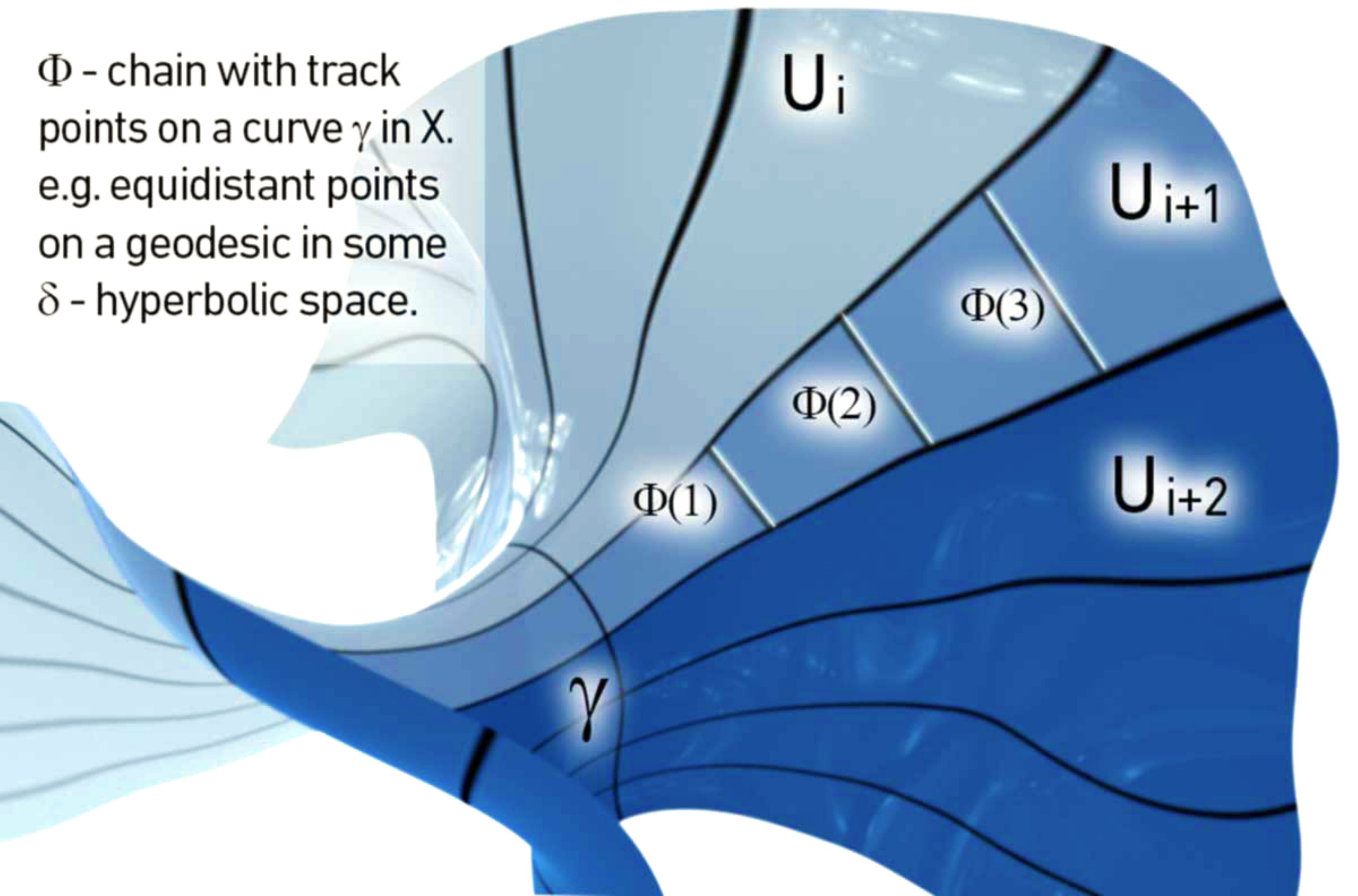}
  \caption{$\Phi$-Chains.}\label{fig:phi}
\end{figure}

\begin{definition}[$\Phi$-Chains]
\label{def:phi}
For a  monotonically increasing function $\Phi:[0,\infty)\to(0,\infty)$ with $\Phi_0:=\Phi(0)>0$ and $\Phi(d)\overset{d\to\infty}{\longrightarrow}\infty$, a \textbf{\textit{Φ}-chain} on a proper geodesic metric space $X$ is a finite or infinite sequence $(U_i)$ of open subsets of $X$ with $U_i\supset U_{i+1}$ together with a sequence of \textbf{track points} $(x_i)$ such that
\begin{enumerate}
\item $\Phi_0\le d(x_i,x_{i+1})\le 3\Phi_0$,
\item $x_i \in \del U_i$,
\item $d(x,\del U_{i\pm1})\ge\Phi(d(x,x_i))$, for every $x\in \del U_i$
\end{enumerate}
for every $i$ where applicable.\footnote{For notational convenience this is slightly different from Ancona's version in \cite[définitions V.5.1]{Anc90}, but essentially the same.} See \autoref{fig:phi}.
\end{definition}

Note that a $\Phi$-chain traversed backwards, i.e., with sets $\cdots\supset X\setminus \overline{U_i}\supset X\setminus \overline {U_{i-1}}\supset\cdots$, is again a $\Phi$-chain with the same track points.
For infinite $\Phi$-chains with arbitrarily large positive indices, $\bigcap_i U_i$ is necessarily empty.

The existence of infinite $\Phi$-chains can be considered as a partial hyperbolicity property of the underlying space. It is easy to see that neither Euclidean space nor asymptotically flat spaces admit any infinite $\Phi$-chains.

A first non-trivial example can be created as follows: with the coordinates $(x,y) \in \R \times \R^{n-1}$ we consider the metric $(1+|y|^2)^2 \cdot g_{\R} +g_\eucl$ on $\R \times \R^{n-1}$. Then the half-spaces $U_i:=(i,\infty) \times \R^{n-1}$ form a $\Phi$-chain with track points $x_i=(i,0)$ for $\Phi(t):=1+t^2$. But as in the Euclidean case, the half-spaces  $U_i[k]:=\R^{k} \times (i,\infty) \times \R^{n-k-1}$, for $1 \le k \le n-1$, do not make up a $\Phi$-chain.

The hyperbolic $n$-space $\Y^n$ carries $\Phi$-chains in all directions and for the same $\Phi$.
Namely, when we think of $\Y^n$ as the upper half-space and $B_{1/2^i}(0)$ is the Euclidean ball of radius $1/2^i$, then $U_i=B_{1/2^i}(0)$, $i \ge 1$, form a $\Phi$-chain for $\Phi(t)=\alpha\cdot t+\beta$ for suitable $\alpha, \beta > 0$, independent even of $n$.  Due to the homogeneity of $\Y^n$ this also gives a $\Phi$-chain along any hyperbolic geodesic $\gamma$, i.e., with track points on $\gamma$.

This ubiquity of $\Phi$-chains, which we see already from considering hyperbolic half-spaces relative to geodesics, extends to arbitrary non-homogenous Gromov hyperbolic spaces.

\begin{theorem}[$\Phi$-Chains on Hyperbolic Spaces]\cite[Section 8]{BHK01}\label{thm:phihyp}
On a proper geodesic $\delta$-hyperbolic space $X$, let $\gamma: [0,4\delta k] \to  X$ be a geodesic with $k\in\Z^+\cup\{\infty\}$. Set $a=\gamma(0)\in X$ and $b=\gamma(4k\delta)\in\overline{X}^{\mathrm{G}}$.
Then the sets
\[
    U_i := \{x\in X\:|\:(x|b)_a> 4i\delta\}
\]

form a $\Phi_\delta$-chain with track points
$x_i=\gamma(4i\delta)$ for $\Phi_\delta(t)=\alpha\cdot t+\beta$, with constants $\alpha=\alpha(\delta)>0$ and $\beta=\beta(\delta)>0$.

\end{theorem}

Note that for a geodesic ray $o=a\leadsto b=\xi\in\delG X$, the sets $U_i=\W_{4i\delta}^o(\xi)\cap X$ are the restriction of our usual neighborhood basis for $\xi$ to $X$.

\subsection{Hyperbolic Unfoldings}\label{sec:hypunf}
A large source for examples of Gromov hyperbolic spaces are \emph{hyperbolic unfoldings} of \emph{uniform spaces}.
The basic idea is that a bounded subspace of a well-controlled space (e.g., a domain in Euclidean space or a compact manifold; the regular part of a singular area-minimizing hypersurface) may degenerate towards a complicated boundary, but on the scale of the distance to this boundary, it might look less obscure. In this situation, one can conformally deform with the inverse of the distance to the boundary to push all difficulties towards infinity. Provided the initial space is \emph{uniform}, this hyperbolic unfolding yields indeed a Gromov hyperbolic space.

In this section, we will sketch the general ideas of this scheme, indicate how Schrödinger operators can be transported to a hyperbolic unfolding and refer to implementations in special cases.

\subsubsection{Uniform Spaces}\label{sec:unifspaces}
Let $(Y, d, \mu)$ be an \emph{incomplete} locally compact metric measure space with metric boundary $\del Y=\overline{Y}\setminus Y$, where $\overline{Y}$ denotes the completion of $Y$. We assume there is a \emph{generalized distance function} $\dd: Y\to(0,\infty)$ on $Y$ which is Lipschitz continuous and vanishes towards $\del Y$. Note that this implies $\dd\le L\,\dist(\cdot,\del Y)$ where $L$ is the Lipschitz constant, and $\dd$ might just be the actual distance to the boundary, but we will later impose further conditions on $\dd$ which could necessitate larger deviations between these functions. Essentially, balls of radius $\sim\dd(y)$ centered at any $y\in Y$ should look uniformly nice.

Such a space is called ($c$-)\emph{uniform}, if each two points can be connected with a $c$-uniform curve, for a global constant $c\ge1$:
\begin{definition}[Uniform Curves]
A rectifiable curve $\gamma: [a,b] \to Y$ is called \textbf{$\bm c$-uniform} with respect to $\dd$  (and $d$), if the following conditions are satisfied:
\begin{itemize}
  \item $\length(\gamma)\le c\, d(x,y)$ ($\gamma$ is a \emph{quasi-geodesic}).
  \item $\min\{\length(\gamma|_{[a,t]}),\length(\gamma|_{[t,b]})\} \le c \, \dd(\gamma(t))$ for any $t\in[a,b]$ (\emph{twisted double-cones}).
\end{itemize}
\end{definition}

A perk of uniform spaces is that they are Gromov hyperbolic in the \emph{\textbf{quasi-hyperbolic metric}}
\[
    k(x,y):=\inf\left\{\int_\gamma1/\dd\:\Bigg|\:\text{$\gamma:x\leadsto y$ rectifiable}\right\}\quad \text{for $x,y\in Y$.}
\]
In the case $\dd=\dist(\cdot,\del M)$, this was first proven by Gehring and Osgood for uniform domains in Euclidean space \cite{GO79} and later generalized to locally compact, rectifiably connected incomplete metric spaces by Bonk, Heinonen and Koskela \cite[Theorem~3.6]{BHK01}. For area-minimizing hypersurfaces with the $\S$-distance $\delta_\ap$, Gromov hyperbolicity of the quasi-hyperbolic metric is proven in \cite[Section~3.2]{Loh18}. This proof uses only axiomatically stated properties of $\delta_\ap$ that are satisfied by our generalized distance function, and inspection of the proof shows that indeed the following is true:

\begin{theorem}[hyperbolization of uniform spaces]\label{thm:grh} If an incomplete locally compact metric space $(Y,d)$ with $L$-Lipschitz generalized distance function $\dd$ is $c$-uniform with respect to $\dd$, it is $\delta$-hyperbolic in the quasi-hyperbolic metric, for some $\delta=\delta(L,c)\ge0$. Moreover, if $(Y,d)$ is bounded\footnote{For unbounded uniform spaces, there is a version of this result involving the one-point compactification \cite[Theorem~3.17]{Loh18}, but we will only be concerned with bounded uniform spaces.}, the metric boundary $\del(Y,d)$ is naturally quasi-symmetrically equivalent\footnote{A homeomorphism $f:(X,d)\to(Y,d')$ between quasi-metric spaces is a \emph{quasi-symmetry}, if there is a homeomorphism $\eta:[0,\infty)\to[0,\infty)$ such that
\[
    \frac{d'(f(x),f(y))}{d'(f(x),f(z))}\le\eta\left(\frac{d(x,y)}{d(x,z)}\right)\quad\text{for distinct points $x,y,z\in X$.}
\]} to the Gromov boundary $\delG(Y,k)$ equipped with a canonical quasi-metric.
\end{theorem}

\begin{remark}\label{rem:smooth}
Even the distance function on Euclidean domains is not everywhere differentiable, and generally in the case of smooth manifolds, it is convenient to have a version of $\dd$ that is smooth, in order to obtain an elliptic operator on the hyperbolization that acts on the usual Sobolev spaces. Such a smooth function $\tdd$, deviating from $\dd$ at most by a constant factor and with prescribed asymptotics of the derivatives, can be constructed by averaging over Whitney cubes, cf.~\cite[VI.2.1, p.\,171]{Ste70}. In the general setup developed so far, this is no more needed, since we can work with appropriately deformed function spaces on the hyperbolization. Note however that in the manifold case, the usage of Whitney smoothings permits to obtain the fundamental results from sections \ref{sec:dirichlet}--\ref{sec:weakcoercivity} with significantly less effort.
\end{remark}

\subsubsection{Schrödinger Forms on Uniform Spaces}\label{sec:operatoronunifspaces}
To apply analytic results on Gromov hyperbolic spaces to a bounded uniform space $Y$, we start with a Schrödinger form $\E=\E^0+V$ which is defined on $\F\subset L^2(Y, \dd^{-2}\mu)$. We expect $\E^0$ to be strongly regular and to satisfy the scale-invariant bounded geometry conditions from \autoref{sec:bg} on balls of radius $\sim\dd$. $V \dd^2$ should be essentially bounded. Then with $N$ denoting the doubling exponent of $(Y,d,\mu)$,
\[
    \E'(u,v):=\E(\dd^{-\frac{N-2}2}u,\dd^{-\frac{N-2}2}v)=\E^0(\dd^{-\frac{N-2}2}u,\dd^{-\frac{N-2}2}v)+\int_Xu\,v\,\dd^2\cdot V\,\d\mu'
\]
defines a Schrödinger form with domain $\F':=\dd^{\frac{N-2}2}\F\subset L^2(Y,\mu':=\dd^{-N}\mu)$ on the complete metric measure space $(Y, k, \mu')$, in particular it satisfies the bounded geometry conditions from \autoref{sec:bg}. $\E'$ is coercive provided $\E$ satisfies a \emph{Hardy inequality}
\[
    \E(u,u)\ge C\int_Yu^2\dd^{-2}\d\mu\quad\text{for every $u\in\F$}
\]
for some $C>0$. This is also called a \emph{strong barrier condition}.

We can bijectively transfer harmonic functions between these two spaces via the mapping $u\mapsto \dd^{\frac{N-2}2}u$, which will allow us to easily carry over results in \autoref{sec:applsing}.

\subsubsection{Examples}
The following spaces and operators satisfy all requirements in the preceding two sections.
\begin{itemize}
\item Uniformity and strong barrier are already interesting conditions for the Laplacian on bounded Euclidean domains with the usual distance to the boundary. They are satisfied on smoothly bounded or Lipschitz domains in Euclidean space \cite{Anc86, Aik12} or more generally domains in spaces of bounded geometry in the sense of \autoref{sec:bg}.
\item Concerning the validity of Hardy inequalities for the Laplacian on metric spaces, the general idea is that the complement of the domain of interest has to be either sufficiently large or sufficiently small. The former case derives from work of Lewis \cite{Lew88}, see \cite{BMS01} for the case of $p$-harmonic functions on metric spaces, while the latter is based on work by Aikawa, see \cite{Leh17} and references therein.
\item The case of Schrödinger operators on the regular part of area-minimizing hypersurfaces is investigated in \cite{Loh18, Loh19, Loh21}. They are uniform with respect to the $\S$-distance $\delta_\ap$. This generalized distance function incorporates information about the second fundamental form to ensure validity of the bounded geometry conditions for the Laplacian. An example of a Schrödinger operator satisfying a Hardy inequality is the \emph{conformal Laplacian} $\L=-\Delta+\frac{n-2}{4(n-1)}\scal$, which admits a strong barrier, given the ambient manifold has nonnegative scalar curvature. As conformal deformation with the first eigenfunction of this operator produces a metric of positive scalar curvature, the Martin theory of such operators is a valuable tool in inductive dimensional descent arguments involving such metrics.
\end{itemize}

\section{Identification of Boundaries}\label{sec:ident}
In this part, we present an adaptation of Ancona's potential theory for weakly coercive operators on Gromov hyperbolic manifolds of bounded geometry \cite{Anc87,Anc90} to the setting of Schrödinger forms on Gromov hyperbolic metric measure spaces. The main result is Ancona's \thmnameref{thm:bhi} which will be used to identify the potential theoretic Martin boundary with the Gromov boundary. Ancona's original article \cite{Anc87} was conceived as a generalization of work of Anderson and Schoen on Cartan--Hadamard manifolds \cite{AS85} and still primarily focuses on those while we will directly approach the more general case of Gromov hyperbolic spaces.

Another valuable source are the comprehensive French lecture notes \cite{Anc90}. They point out connections to heat kernels and stochastic processes such as Brownian motion and introduce potential theory on graphs. Our setting permits us to unify the cases of graphs and manifolds, which Ancona considered separately. Regarding the apparent similarities between the arguments in these cases, it is no surprise that our streamlined account requires few substantial modifications. Additionally, we carefully keep track of all involved constants to show that the quantitative results depend only on a set of universal constants, not on the explicit space or operator under consideration. This is useful for blow-up arguments such as in \cite{Loh19,Loh21} where these constant can be shared amongst sequences and limit spaces.

\subsection{Hyperbolicity and Boundary Harnack Inequalities}\label{sec:hy}

Now we additionally invest the hyperbolicity of the underlying space $X$. We employ the property that in Gromov hyperbolic spaces any two points can be connected by well-controlled $\Phi$-chains as in \thmref{thm:phihyp}. We prove in fact a more general result which holds on a single $\Phi$-chain, even if the space carries essentially only this one $\Phi$-chain as in the example $(\R \times \R^{n-1},(1+|y|^2)^2 \cdot g_{\R} +g_\eucl)$. This can be seen as a directed form of hyperbolicity.

The technical main result of this section describes the behavior of Green function along $\Phi$-chains, building on the results from \autoref{sec:glob}. In the presence of sufficiently many $\Phi$-chains, such as in Gromov hyperbolic spaces, we infer boundary Harnack inequalities and employ them to identify the Martin boundary with the Gromov boundary.

Our general assumptions remain the same as in the previous sections: $\E$ is a coercive Schrödinger form with associated operator $\L$ and Green function $G$ on a complete metric measure space $(X,d,\mu)$. The additional assumptions, that is, the presence of a $\Phi$-chain (depending on the function $\Phi$) or even of an underlying hyperbolic geometry (with constant $\delta$ and coming with a universal function $\Phi=\Phi_\delta$), are stated directly in the results.

\subsubsection{Global Behavior: \texorpdfstring{\textit{Φ}}{Phi}-Chains}\label{sec:gb}

The following result, sometimes called \emph{3G-inequality}, describes the major influence of $\Phi$-chains on the behavior of Green functions.

\begin{theorem}[Green functions along $\Phi$-chains]
\label{thm:gphi}
There is a suitable constant $c(\sigma, N, C_P, C_D, k, \epsilon,\Phi) >1$ such that for any $\Phi$-chain with track points $x_1,\dots,x_m$ (as in \autoref{def:phi}), we have for the minimal Green function $G$ the estimate
\begin{align*}
c^{-1}\mu(B_\sigma(x_j))\,G(x_m,x_j)\,G(x_j,x_1)&\le G(x_m,x_1)\\
&\le c\,\mu(B_\sigma(x_j))\,G(x_m,x_j)\,G(x_j,x_1)\,,\quad j=2,...,m-1.
\end{align*}
\end{theorem}

At the heart of the argument we employ the pairing of two at first sight entirely unrelated geometric and analytic properties: the existence of $\Phi$-chains on $X$  and coercivity of $\L$.
The idea is that, on the one hand, $\Phi$-chains allow finding balls of arbitrary large radii in $U_{i-1} \setminus U_{i+1}$ centered in $\p U_i$ within a uniformly upper bounded distance to the track points. On the other hand, the \thmnonameref{qmpglob}{relative maximum principle} shows that on these balls we can improve crude estimates from a previous application of a
Harnack inequality by investing coercivity. This makes the following result the main step in the proof of the Theorem.

\begin{proposition}[growth recovery along $\Phi$-chains]
\label{prop:expg}
For any given $\Phi$-chain with track points $x_1,\dots,x_j$, we have
\begin{equation}
\label{eq:step1}
G(z,x_1)\le c\,\mu(B_\sigma(x_j))\,G_{-\epsilon}(z,x_j)\,G(x_{j+1},x_1)\quad\text{for }z\in\del U_{j+1},
\end{equation}
for some constant $c(\sigma, N, C_P, C_D, k, \epsilon,\Phi) >0$ independent of the length $j$ of the $\Phi$-chain.

\begin{figure}[h]
  \centering
  \includegraphics[width=.8\textwidth]{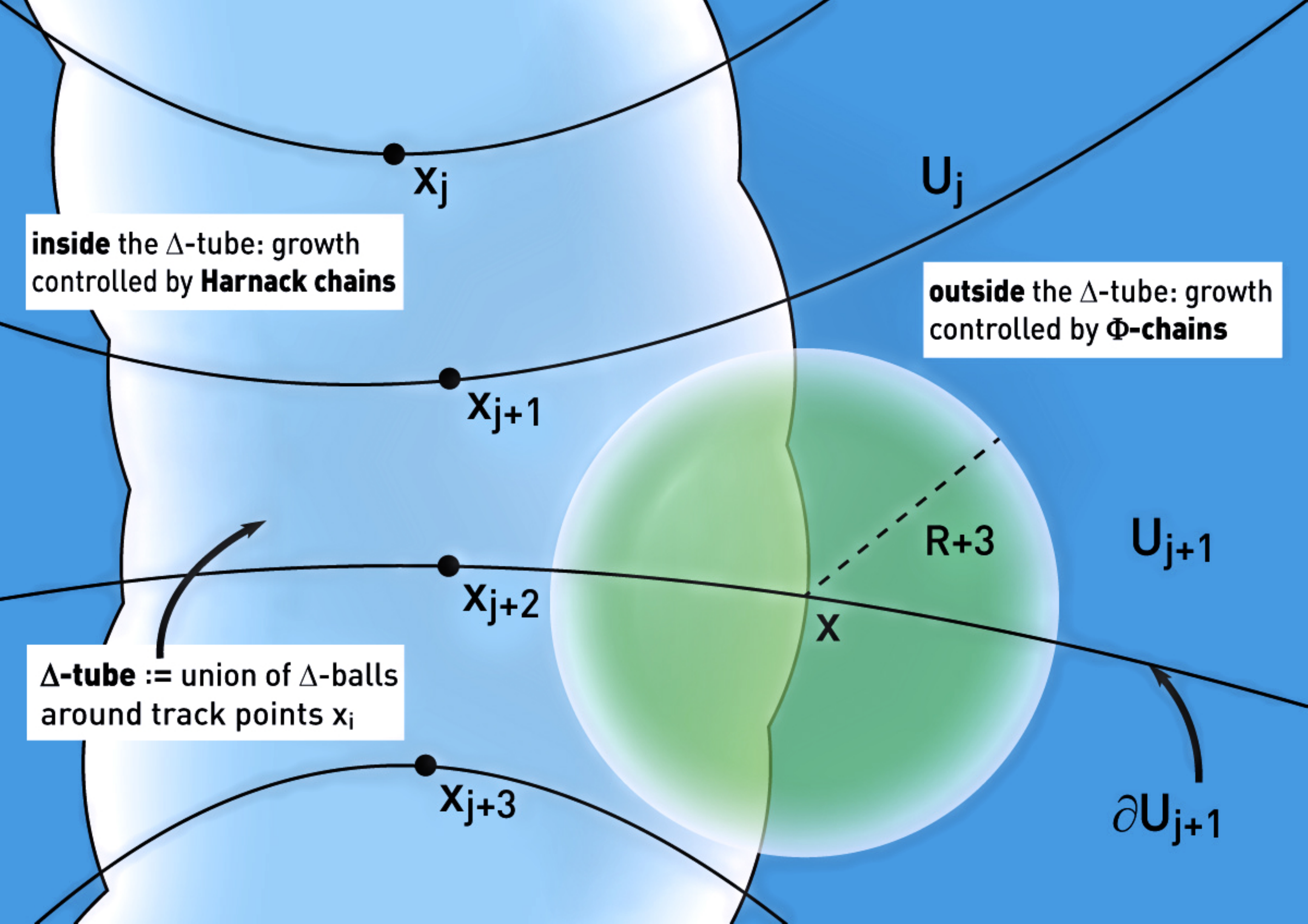}
  \caption{Growth Recovery Mechanism.}
\end{figure}

\begin{proof} The argument is by induction over the length $j$.

For $j=1$, we can use $G_{-\epsilon}\ge G$, inequality \eqref{eq:resineq1}, the \thmnonameref{prop:bgreen}{lower bound for the Green function} $c_1\,\mu(B_\sigma(x_1))\,G(x_2,x_1)\ge1$ and Harnack inequalities to get a first guess for the constant $c=c_1$ and note that $c_1$ depends only on $\sigma$, $N$, $C_P$, $C_D$, $k$, $\epsilon$ and $\Phi$ (via Harnack inequalities).

For the induction step we first assume we have proved the weaker assertion that there is a constant $c_j$ so that \eqref{eq:step1} holds for any $\Phi$-chain of length $j$.  Then we can apply the Harnack inequalities for $\L$ and $\L_{-\epsilon}$ and doubling of $\mu$ to get a constant $c'(\sigma, N, C_P, C_D, k, \epsilon, \Phi)\ge 1$, independent of $j$, such that
\begin{equation}
\label{eq:step1a}
G(z,x_1)\le c'c_j\,\mu(B_\sigma(x_{j+1}))\,G_{-\epsilon}(z,x_{j+1})\,G(x_{j+2},x_1)\quad\text{for }z\in\del U_{j+1}\,.
\end{equation}
By the \thmnameref{gmp} this inequality extends to $z\in \bar U_{j+1}$. Now we invest the coercivity of $\L$ and the properties of the $\Phi$-chains to improve this inequality.\medskip

Towards this end, we first apply the \thmnonameref{qmpglob}{relative maximum principle} to the ($\L$-superharmonic) function $G_{-\epsilon}(\cdot,x_{j+1})$ and its greatest $\L$-harmonic minorant $u$ on some ball $B_R(x)$ which we can represent as $u=\red_{G_{-\epsilon}(\cdot,x_{j+1})}^{\del B_R(x)}$, the reduit always taken with respect to $\L$ (see \autoref{sec:balayage} for properties of the reduit employed here and below). For $R=\ln(1/c')/\ln \eta$ and $B_{R+\sigma}(x) \subset U_{j+1}$, the relative maximum principle yields
\[ u(x)\le \frac1{c'}G_{-\epsilon}(x,x_{j+1})\,.\]

In turn, the definition of a $\Phi$-chain shows that there is some $\Delta(\Phi,R)>0$ such that $B_{R+3}(x)\subset U_{j+1}$, for $x\in\del U_{j+2}$, as soon as $d(x,x_{j+2})\ge\Delta$. Then, we have from  \eqref{eq:step1a}
\begin{align*}
G(x,x_1)&=\red_{G(\cdot,x_1)}^{\del B_R(x)}(x)\le c'c_j \mu(B_\sigma(x_{j+1})) \red_{G_{-\epsilon}(\cdot,x_{j+1})}^{\del B_R(x)}(x)\,G(x_{j+2},x_1)\\
&=c' c_j \mu(B_\sigma(x_{j+1})) u(x)\, G(x_{j+2},x_1) \le  c_j \mu(B_\sigma(x_{j+1})) G_{-\epsilon}(x,x_{j+1})\, G(x_{j+2},x_1).
\end{align*}

On the other hand, we get universal estimates for $x\in\del U_{j+2}$ with $d(x,x_{j+2})<\Delta$: since $d(x,x_{j+2})<\Delta$, there are constants $c'',c'''\ge1$ only depending on $\Delta$, $\Phi$, $H$ and $c_1$, such that $G(x,x_1)\le c''\, G(x_{j+2},x_1)$ by Harnack inequalities and $\mu(B_\sigma(x_{j+1}))G_{-\epsilon}(x,x_{j+1})\ge (c''')^{-1}$ by the \thmnonameref{prop:bgreen}{bounds for the Green function} and Harnack inequalities. The result is
\[ G(x,x_1)\le c'' c''' \,\mu(B_\sigma(x_{j+1}))\,G_{-\epsilon}(x,x_{j+1})\,G(x_{j+2},x_1)\quad\text{for $x\in\del U_{j+2}$ with $d(x,x_{j+2})<\Delta$}\,.\]

Everything combined, we can choose $c=\max\{c_1,c''c'''\}$ and outside a tube of radius $\Delta$ around the track points the constant can be kept in every induction step while on the inside we can always use the universal constant.
\end{proof}
\end{proposition}

\begin{proof}[Proof of {\thmref{thm:gphi}}]
The first inequality is rather easy:
For $x\in\del B_\sigma(x_j)$ we have
\[ \mu(B_\sigma(x_j))G(x,x_j)\,G(x_j,x_1)\le c_1 G(x_j,x_1)\le c_1 H\,G(x,x_1) \]
by \thmref{prop:bgreen} and the Harnack inequality. Since the left-hand side is an $\L$-potential and the right-hand side is $\L$-superharmonic, this inequality extends to $X\setminus B_\sigma(x_j)$ and in particular to $x_m$ by the \thmnameref{gmp}.

For the second inequality, we use repeatedly \thmref{prop:expg} and the \thmnameref{resolvent}:
\begin{align}
G(x_m,x_1)&=\red_{G(\cdot,x_1)}^{\del U_j}(x_m) &&|\ x_1\notin U_j\nonumber\\
 &\overset{\eqref{eq:step1}}{\le}c\,\mu(B_\sigma(x_j))\,\red_{G_{-\epsilon}(\cdot,x_j)}^{\del U_j}(x_m)\, G(x_{j+1},x_1)\nonumber\\
 &=c\,\mu(B_\sigma(x_j))\,\red_{G(\cdot,x_j)+\epsilon G(G_{-\epsilon}(\cdot,x_j))}^{\del U_j}(x_m)\,G(x_{j+1},x_1)  &&|\text{ res.eq.}\nonumber\\
 &\le c\,\mu(B_\sigma(x_j))\left(G(x_m,x_j)+\epsilon \int_X \red_{G(\cdot,z)}^{\del U_j}(x_m)\,G_{-\epsilon}(z,x_j)\,\d\mu(z) \right)G(x_{j+1},x_1)\label{eq:long}
\end{align}
At this point, we can again employ the first step \eqref{eq:step1}, but now for the reversed $\Phi$-chain $x_m,\dots,x_1$ with $X\setminus U_m,\dots,X\setminus U_1$, namely
\[ G(x_m,z)\le c\,\mu(B_\sigma(x_{j+1}))\,G_{-\epsilon}(x_{j+2},z)\,G(x_m,x_{j+1})\quad\text{for }z\in X\setminus U_{j+1}\,.\]
This holds on all of $X\setminus U_{j+1}$ by the \thmnameref{gmp}. Since $x_j\in X\setminus U_{j+1}$, this can be directly applied to $G(x_m,x_j)$.
For the second summand in \eqref{eq:long}, we have $\red_{G(\cdot,z)}^{\del U_j}(x_m)=\red_{G(x_m,\cdot)}^{\del U_j}(z)\le G(x_m,z)$ for $z\in \del U_j\subset X\setminus U_{j+1}$ (see \autoref{sec:balayage}), but the upper bound $\red_{G(x_m,\cdot)}^{\del U_j}(z)\le c\,\mu(B_\sigma(x_{j+1}))\,G_{-\epsilon}(x_{j+2},z)\,G(x_m,x_{j+1})$ is valid for \emph{all} $z\in X$ by definition of the reduit since the right-hand side is positive and $\L$-superharmonic in $z$.
Thus,
\begin{align*}
G(x_m,x_1) &\le c^2\,\mu(B_\sigma(x_{j+1}))\,G(x_m,x_{j+1})\,G(x_{j+1},x_1)\\
&\quad\cdot\left[\mu(B_\sigma(x_j))\left(G_{-\epsilon}(x_{j+2},x_j)+\epsilon\int_X G_{-\epsilon}(x_{j+2},z)\,G_{-\epsilon}(z,x_j)\,\d\mu(z)\right)\right].
\end{align*}
The large square bracket is universally bounded from above by \thmref{prop:bgreen}, the Harnack inequality, and the inequalities \eqref{eq:resineq1} and \eqref{eq:resineq2} following from the \thmnameref{resolvent} for $t=\frac32\epsilon$.
\end{proof}

Now we assume that $X$ is a $\delta$-hyperbolic space, then we can choose $\Phi=\Phi_\delta$ and recall that there are $\Phi$-chains along geodesics in $X$. Since $\Phi_\delta$ is
determined from $\delta$, the $\Phi$-dependence of the estimates now reduces to a $\delta$-dependence.

\begin{corollary}[Green function along hyperbolic geodesics]
\label{thm:gfalg} If $X$ is $\delta$-hyperbolic, let $x,y,z\in X$ such that $y$ lies on geodesic connecting $x$ and $z$ with $d(x,y),d(y,z)>22\delta$.
Then there is a constant $c(\sigma, N, C_P, C_D, k, \epsilon,\delta) >1$ such that
\[
    c^{-1}\,\mu(B_\sigma(y))G(x,y)\,G(y,z)\le G(x,z)\le c\,\mu(B_\sigma(y))\,G(x,y)\,G(y,z)\,.
\]
\end{corollary}

\begin{remark}
This estimate for the Green function can be interpreted stochastically and algebraically:
\begin{itemize}
\item Stochastically, the Green function $G(x,y)$ is a density for the expected number of times an $\L$-Brownian motion starting at $y$ reaches $x$, see \cite{Pin95} or \cite{Anc90}. Now the estimate for the Green function above states that on a hyperbolic geodesic $x\leadsto y\leadsto z$, up to a constant multiple there are as many Brownian particles travelling directly from $x$ to $z$ as there are particles travelling from $x$ to $y$ and then to $z$. This is in line with the geometric \enquote{valley} interpretation of Gromov hyperbolic spaces, as seen in \thmref{prop:geostab}.
\item Algebraically, we can examine the function
\[
    d_G(x,y):=-\ln\left(G(x,y)\sqrt{\mu(B_\sigma(x))\mu(B_\sigma(y))}\right)\quad\text{for $x,y\in X$.}
\]
For this \emph{Green metric}, the estimate above says
\[
    d_G(x,z)=d_G(x,y)+d_G(y,z)\pm\ln c
\]
along a geodesic $x\leadsto y\leadsto z$, if the points are sufficiently far apart. Note that the estimate \enquote{$\le$} holds not only along a geodesic, but for general points $x$, $y$, $z$ with large mutual distance -- this is a rough triangle inequality for $d_G$. The other direction suggests that the large-scale geodesic structure of $d_G$ is comparable to that of $d$. In the context of hyperbolic groups, this has been explored further in \cite{BB07,BHM11}, where the Green metric turns out to be Gromov hyperbolic and quasi-isometric to the word metric for non-amenable groups and certain operators.
\end{itemize}
\end{remark}

\subsubsection{Boundary Harnack Inequality}
\label{sec:bhp}
To formulate a boundary Harnack inequality near points on the Gromov boundary of a $\delta$-hyperbolic space, we use the following Euclidean version as a blueprint:
\begin{theorem}[boundary Harnack inequality on a disc]\label{thm:bhid}\cite{Kem72}
There exist constants $A$, $C>1$ such that for any point $\xi\in \p B_1(0)\subset\R^2$ and $0<R<1$ the following is true: for any two harmonic functions $u$, $v>0$ (with respect to the Laplacian) on $B_{A \cdot R}(\xi) \cap B_1(0)$ that vanish along $B_{A \cdot R}(\xi) \cap \p B_1(0)$,
\[
\frac{u(x)}{v(x)} \le C   \frac{u(y)}{v(y)} \text{ for all } x,\, y \in B_R(\xi) \cap B_1(0)\,.
\]
\end{theorem}
Note that in the non-boundary version of the Harnack inequality, the appearance of another function $v$ is obscured by the fact that the constant function 1 is (nearly) harmonic.

As a replacement for balls in the classical version of the boundary Harnack inequality, we need some characterization of neighborhoods of a point at infinity. This is made precise by the notion of $\Phi$-neighborhood bases.
\begin{definition}[$\Phi$-neighborhood basis]
Two non-empty open subsets $V\supset W$ of $X$ are called \textbf{$\bm{\Phi}$-neighborhoods with hub $h\in X$}, if $\overline{W}\subset V$, $B_{\Phi_0}(h) \subset V\setminus\overline{W}$ and any two points $p \in \del V$ and $q \in \del W$ can be joined by
a $\Phi$-chain that has $h$ as a track point. We call an infinite family of non-empty open sets $\N_i\subset X$, $i=1,2,3,\dots$, with $\bigcap_i\N_i=\ep$ a \textbf{$\bm{\Phi}$-neighborhood basis}, if $\N_i$ and $\N_{i+1}$ are $\Phi$-neighborhoods with hub $o_i$, for every $i$.
\end{definition}
Just as metric balls are, besides their role as neighborhood bases, the basic playground for Harnack inequalities, these $\Phi$-neighborhoods are their counterpart in \emph{boundary} Harnack inequalities.

In $\delta$-hyperbolic spaces, every point in the Gromov boundary has a canonical neighborhood basis that is also a $\Phi_\delta$-neighborhood basis. Namely, as shown in \cite[Proposition 8.10]{BHK01}\footnote{They define the neighborhood basis with references to distances from geodesics, but the two concepts can easily be translated into each other using \thmref{lem:productasdistance}.}, we have:

\begin{lemma}[$\Phi_\delta$-neighborhood basis] \label{lemma:den} If $X$ is $\delta$-hyperbolic and $o\in X$ a basepoint, there is a constant $c_\delta>0$ only depending on $\delta$ such that for any $\xi \in \delG X \subset\overline{X}^{\mathrm{G}}$, the open sets
\[
    \N^\delta_i(\xi):=\W_{c_\delta i}^o(\xi)\cap X=\{x\in X\:|\:(x|\xi)_o>c_\delta i\}\quad\text{for $i=1,2,\dots$}
\]
are a $\Phi_\delta$-neighborhood basis. Recall that their closures $\overline{\N^\delta_i(\xi)}=\overline{W_{c_\delta i}^o(\xi)}\subset \overline{X}^{\mathrm{G}}$ in the Gromov compactification form a neighborhood basis of $\xi\in\delG X$.
\end{lemma}

For later applications on uniform domains, it is not always possible to find $\L$-harmonic functions that vanish at infinity, because there are situations where even minimal Green functions might diverge. Hence, we introduce a more robust, operator-aware notion of vanishing which can be thought of as a minimal growth condition. This will be further explained in \thmref{prop:LvanM} below.

\begin{defprop}[$\L$-vanishing]\label{Lvan}
We say that a positive $\L$-superharmonic function $u$ \emph{\textbf{$\bm{\L}$-vanishes}} on an open set $V\subset X$ if one of the following equivalent conditions is satisfied:
\begin{enumerate}
\item There is a positive $\L$-superharmonic function $w$, such that $u/w\to0$ at infinity, i.e., for every $\epsilon>0$ there is a compact set $K\subset X$ with $u/w<\epsilon$ on $V\setminus K$.
\item There is an $\L$-potential $p$ such that $p\ge u$ on $V$.
\item  The reduit $\red_u^V$ is an $\L$-potential (see \autoref{sec:balayage}).\footnote{\cite{Anc87} uses the first definition, but it is easier to employ the last.}
\end{enumerate}
\begin{proof}
$\mathrm{(i)}\Rightarrow\mathrm{(iii)}$
Assume there is a positive $\L$-harmonic function $h$, such that $\red_u^V\ge h>0$ on $X$. For some fixed $\epsilon>0$, choose a compact set $K$ with $u<\epsilon w$ on $V\setminus K$. By the properties of the reduit, we even have
\[\epsilon w\ge\red_u^{V\setminus K}\ge\red_u^V-\red_u^{V\cap K}\ge h-\red_u^{V\cap K}
\]
on all of $X$. Now $\red_u^{V\cap K}$ is an $\L$-potential since $V\cap K$ is relatively compact in $X$, $\epsilon w-h$ is $\L$-superharmonic and $\epsilon w\ge u=\red_u^V\ge h$ on $\del(V\setminus K)=\del(X\setminus(V\setminus K))\subset X$, thus we can apply the \thmnameref{gmp} to see $\epsilon w\ge h$ on all of $X$. Since $\epsilon$ was arbitrary, $h=0$.
\medskip

$\mathrm{(iii)}\Rightarrow\mathrm{(ii)}$
Choose $p=\red_u^V$.
\medskip

$\mathrm{(ii)}\Rightarrow\mathrm{(i)}$
It suffices to show that an $\L$-potential $p$ satisfies the condition everywhere. Towards this end, consider the functions $\red_p^{X\setminus \overline{B_R(o)}}$ on balls around an arbitrary basepoint $o\in X$. They converge to zero for $R\to\infty$ because the limit is $\L$-harmonic and $\le p$. Let $(x_j)$ be a countable dense set in $X$. We may choose a sequence $(R_i)$ such that $\red_p^{X\setminus \overline{B_{R_i}(o)}}(x_j)\le 2^{-i}$ for all $j\le i$. Then the function $w=\sum_i\red_p^{X\setminus \overline{B_{R_i}(o)}}$ is finite on a dense set, $\L$-superharmonic and $p/w\le 1/i$ outside $B_{R_i}(o)$.
\end{proof}
\end{defprop}

Note that all $\L$-potentials such as the minimal Green function are $\L$-vanishing on $V=X$ and hence on all open sets because the property of $\L$-vanishing is conserved on subsets as can be easily seen using condition (i).

On bounded sets, $\L$-vanishing \emph{at infinity} is trivially true for any $\L$-superharmonic function. The concept is also useless for unbounded sets shrinking too quickly when approaching infinity. On Gromov hyperbolic spaces, it becomes significant for sets $V=W \cap X$, where $W\subset \overline{X}^\mathrm{G}$ is open with non-empty intersection $W \cap \delG X$. For more elaborate criteria in the context of Martin theory see \thmref{prop:LvanM}.
\medskip

On $\Phi$-neighborhoods we can now formulate the following central result.
\begin{theorem}[boundary Harnack inequality]
\label{thm:bhi}
Let $V\supset W$ be $\Phi$-neighborhoods with hub $h$ and $u$, $v$ two positive $\L$-superharmonic functions that are $\L$-harmonic and $\L$-vanishing on $V$, then there is a constant $H_B=H_B(\sigma, N, C_P, C_D, k, \epsilon, \Phi)$ such that
\[
    \frac{u(x)}{u(y)}\le H_B \frac{v(x)}{v(y)}\quad\text{for any $x,y\in W$}\,.
\]
\begin{proof}
By \thmref{Lvan}, the reduit $\red_u^{V}$ is an $\L$-potential and therefore admits a representation as
\[ \red_u^{V}=\int_{\del V} G(\cdot,z)\,\d\nu(z) \]
for some Radon measure $\nu$. On $\overline{W}$, $\red_u^{V}$ agrees with $u$, and therefore we have
\begin{align*}
  u(x)&=\int_{\del V} G(x,z)\,\d\nu(z)\\
  &\le c\,\mu(B_\sigma(h))\int_{\del V} G(x,h)\,G(h,z)\,\d\nu(z)\\
  &= c\,\mu(B_\sigma(h))\,G(x,h)\,u(h)\text{ for $x\in\del W$}
\end{align*}
using the assumption that $x\in\del W$ and $z\in\del V$ can be connected by a $\Phi$-chain through $h$ and \thmref{thm:gphi}. The other inequality from there gives
\[
  v(x)\ge c^{-1}\,\mu(B_\sigma(h))\,G(x,h)\,v(h)\quad\text{for $x\in\del W$}
\]
and both inequalities extend to $W$ by the \thmnameref{gmp} because $\red_u^{V}$ and $G(\cdot,h)$ respectively are $\L$-potentials. We can combine them to obtain
\[ \frac{u(x)}{v(x)}\le c^2\frac{u(h)}{v(h)}\quad\text{for $x\in W$}\,.\]
Interchanging the roles of $u$ and $v$ yields the result with $H_B=c^4$.
\end{proof}
\end{theorem}

In most applications, $u$ and $v$ are either globally $\L$-harmonic functions or minimal Green functions with pole outside $V$.
\medskip

In the case of a $\delta$-hyperbolic space, the size of the smaller neighborhood can be explicitly quantified using \thmref{lemma:den}. We get

\begin{corollary}[hyperbolic boundary Harnack inequality]
\label{cor:hbhi} If $X$ is $\delta$-hyperbolic, there is some positive constant $H_B(\sigma, N, C_P, C_D, k, \epsilon, \delta)>1$ such that two positive $\L$-superharmonic functions $u$, $v$ that are $\L$-harmonic and $\L$-vanishing on a $\Phi_\delta$-neighborhood $\N^\delta_i(\xi)$ of $\xi\in\delG X$ satisfy
\[ \frac{u(x)}{u(y)}\le H_B \frac{v(x)}{v(y)}\quad\text{for any }x,y\in \N^\delta_{i+1}(\xi)\,. \]
\end{corollary}

%%%%%%%%%%%%%%%%%%%%%
%% MARTIN BOUNDARY %%
%%%%%%%%%%%%%%%%%%%%%
\subsubsection{The Hyperbolic Martin Boundary}
We now turn our attention towards the \emph{Martin boundary} as introduced in \autoref{sec:martin}.
In the situation at hand, $\Phi$-neighborhood bases are essentially neighborhood bases of minimal Martin boundary points.
\begin{theorem}[characterization of minimal Martin points]\label{cmmp}
Assume $o\in X$ is a basepoint.
Let $(\N_i)$ be a $\Phi$-neighborhood basis with hubs $h_i$ and assume $o\notin\overline{\N_1}$.
Denoting the interior of the closure of $\N_i\subset X \subset \overline{X}^\mathrm{M}$ in the Martin compactification $\overline{X}^\mathrm{M}$ of $X$ by $\widetilde{\N}_i$, there is exactly one Martin boundary point $\zeta$ in $\bigcap\widetilde{\N}_i$. The resulting Martin function $K_\zeta$ is characterized as the only positive $\L$-harmonic function $\L$-vanishing on every $X\setminus\N_i$, up to scalar multiples. In particular, this Martin point $\zeta$ is minimal.
\begin{proof}
By the Harnack inequalities, the sequence $K_{h_i}=\frac{G(\cdot, h_i)}{G(o,h_i)}$ has a subsequence locally uniformly converging to some $\L$-harmonic function $K_\zeta$ representing a Martin boundary point $\zeta$.
$K_\zeta$ is $\L$-vanishing on every $X\setminus\N_i$ because by the \thmnameref{thm:bhi}, every $K_{h_j}$ for $j\ge i$, and hence every limit, is upper bounded by the $\L$-potential $H_B K_{h_i}$ on $X\setminus\N_i$.

Now assume there is another positive $\L$-harmonic function $u$ that is $\L$-vanishing on $X\setminus\overline{\N}_i$ for every $i$, w.l.o.g.\ $u(o)=1$. Applying the \thmnameref{thm:bhi} we see $H_B^{-1}K_\zeta\le u \le H_B K_\zeta$ on $X$.
Thus, $\eta:=\inf u/K_\zeta\le1$ is positive. By the strong minimum principle following from the \thmnameref{harnackineq}, the $\L$-harmonic function $u-\eta K_\zeta\ge0$ has to be positive everywhere, else it would be identical zero.
In the former case we can again apply the \thmnameref{thm:bhi} to $K_\zeta$ and $u-\eta K_\zeta$ to get $(\eta+(1-\eta)H_B^{-1})K_\zeta\le u$ which contradicts the definition of $\eta$, unless $\eta=1$ and $u=K_\zeta$.
\end{proof}
\end{theorem}

This gives the following characterization of the Martin boundary in case we have enough $\Phi$-neighborhoods:
\begin{corollary}[identifying the Martin boundary]\label{imb}
If in a given compactification $\overline X$ of $X$ (i.e., $\overline X$ is compact and $X\subset\overline X$ open and dense) every boundary point admits a neighborhood basis of the form $\overline\N_i\subset \overline X$ for some $\Phi$-neighborhood basis $(\N_i)$, it is canonically homeomorphic to the Martin compactification $\overline X^\mathrm{M}$.
\begin{proof}
\thmref{cmmp} yields an injective map from $\overline X$ to $\overline X^\mathrm{M}$. It is continuous because for every sequence $(y_i)$ in $\overline X$ converging to $\zeta\in\overline X\setminus X$ the corresponding Martin functions $K_{y_i}$ converge to the unique Martin function that is $\L$-vanishing on all $X\setminus\N_i$ for some $\Phi$-neighborhood basis $(\N_i)$ of $\zeta$, that is $K_\zeta$. Thus, by elementary properties of compactifications \cite[§38]{Mun00}, it is already a homeomorphism. Note that in particular all Martin boundary points are minimal.
\end{proof}
\end{corollary}

From this and \thmref{lemma:den} we get the following principal potential theoretic result on Gromov hyperbolic spaces.
\begin{corollary}[Gromov boundary and Martin boundary] \label{cor:gmr}   Assume that $X$ is Gromov hyperbolic. Then the Gromov and Martin boundaries of $X$ are canonically homeomorphic and every Martin boundary point is already minimal,
\[ \delG X\cong \delM (X,\L)\cong  \delM^0 (X,\L) .\]
This means that a positive function $u>0$ on $X$ is $\L$-harmonic if and only if there is a (unique) Radon measure $\mu_u$ on $\delG X$
such that
\[
    u(x)  =\int_{\delG X} K_\zeta(x) \, \d \mu_u(\zeta).
\]
\end{corollary}

Now that we have the right notion for a potential theoretic boundary at infinity, we can update the notion of $\L$-vanishing. Note that the following characterizations are slightly different from the formulation in \thmref{Lvan} because there, without explicit references to a boundary, it was only possible to refer to $\L$-vanishing \emph{at infinity of open sets $V$ in $X$}, i.e., on the Martin boundary points in $\overline{V}\cap\delM X\subset\overline{X}^\mathrm{M}$.
\begin{proposition}[$\L$-vanishing and Martin boundary]\label{prop:LvanM}
Assume every Martin boundary point has a $\Phi$-neighborhood basis, e.g., $X$ is Gromov hyperbolic.
For an open subset $\varXi\subset \delM X$ of the Martin boundary and a positive $\L$-harmonic function $u$ on $X$ the following are equivalent:
\begin{enumerate}
\item $u$ is $\L$-vanishing on any open set $V\subset X$ with $\overline V\cap\delM X\subset\varXi$ in the Martin compactification.
\item On any open set $V\subset X$ with $\overline V\cap\delM X\subset\varXi$ in the Martin compactification, the following property holds: each positive $\L$-harmonic function $v$ on $V$ with $v\ge u$ on $\del V\cap X$ satisfies $v\ge u$ on $V$.
\item The Martin measure $\mu_u$ associated to $u$ is supported outside $\varXi$, i.e., $\mu_u(\varXi)=0$.
\end{enumerate}
\begin{proof}
$\mathrm{(i)}\Rightarrow\mathrm{(ii)}$
If $u$ is $\L$-vanishing on $V$, $\red_u^V$ is an $\L$-potential and $\red_u^V=u$ on $V$. Then the \thmnameref{gmp} gives exactly the desired property.
\medskip

$\mathrm{(iii)}\Rightarrow\mathrm{(i)}$
Each Martin function $K_{\zeta}$ is $\L$-vanishing outside $\zeta$ and, therefore, the Martin integral representing $u$ $\L$-vanishes outside the support of  $\mu_u$.
\medskip

$\mathrm{(ii)}\Rightarrow\mathrm{(iii)}$
Assume $\mu_u(\varXi)\neq0$, then there is a compact $K\subset\delM X$ and an open set $W\subset\overline{X}^\mathrm{M}$ such that $K\subset W\cap\delM X\Subset\varXi$, $V:=W\cap X$ satisfies the condition in (ii), and $\mu_u(K)\neq0$ (since the Radon measure $\mu_u$ is inner regular). Therefore, it is enough to consider the case where $u \equiv u_K:=\int_K K_\zeta\,\d\mu_u(\zeta)$.

We compare $u_K$ with the minimal Green function $G(\cdot,p)$ with pole $p\in X\setminus V$. Recalling the argument of $\mathrm{(iii)}\Rightarrow\mathrm{(i)}$ we know that $u_K$ is $\L$-vanishing on an open neighborhood $N$ of $X\setminus V$. By compactness of $\overline{X\setminus V}$ in the Martin compactification, $X\setminus V$ can be covered by finitely many $\Phi$-neighborhoods contained in $N$ and a compact subset of $X$. Then the \thmnameref{thm:bhi} shows that there is a $C>0$ such that $C\cdot G(\cdot,p)\ge u_K$ on $X\setminus V$ and especially on $\del V$. But then from (ii) it follows that $C\cdot G(\cdot,p)\ge u_K$ on $V$, hence on all of $X$, hence $u_K\equiv 0$ because $G(\cdot,p)$ is an $\L$-potential.
\end{proof}
\end{proposition}

\subsubsection{Sharpness}
The validity of such a simple identification of the Martin boundary with a geometric boundary, which one may expect from a naive guess, actually is a rare exception. If we only slightly violate the hyperbolicity constraint we can get a completely different and rather inscrutable outcome with many non-minimal Martin boundary points.
\begin{example}[Ideal Boundaries of $\mathbb{H}^m \times \mathbb{H}^n$] \label{ex:pro}
The product space of two classical hyperbolic spaces  $\mathbb{H}^m \times \mathbb{H}^n$, $m, n \ge 2$, has bounded geometry, but it is no longer Gromov hyperbolic since it contains flat planes obtained as products of pairs of geodesics in the two factors. In \cite{GW93} we find a thorough discussion of the case of the Laplace operator on $\mathbb{H}^m \times \mathbb{H}^n$. There exist positive functions $s$ with $-\Delta s = \lambda \cdot s$ if and only if $\lambda \le \lambda_0$, where $\lambda_0$ is the principal eigenvalue. For $\lambda<\lambda_0$, the operator $-\Delta - \lambda$ is coercive.
The boundary at infinity (a natural generalization of Gromov boundary) is homeomorphic to $S^{n+m-1}$ \cite[II.8.11(6), p.\,266]{BH99}. In turn, for the minimal Martin boundary of $-\Delta - \lambda$, for $\lambda<\lambda_0$,  we have
\[\delM^0 (\mathbb{H}^m \times \mathbb{H}^n,-\Delta - \lambda) = S^{m-1} \times  S^{n-1} \times I_\lambda\,,\]
where $I_\lambda$ is a closed interval with a natural parametrization which depends on $\lambda$ and degenerates to a single point when $\lambda \ra \lambda_0$ \cite[p.\,21]{GW93}. The full Martin boundary contains two additional pieces,
\begin{equation}\label{part}
\delM (\mathbb{H}^m \times \mathbb{H}^n,-\Delta - \lambda) = (S^{m-1} \times \mathbb{H}^n) \cup (S^{m-1} \times  S^{n-1} \times I_\lambda) \cup (\mathbb{H}^m \times S^{n-1})/\sim\,.
\end{equation}
The gluing maps for $\sim$ are described in \cite[p.\,27]{GW93}. We observe that not only the boundary at infinity does not coincide with $\delM^0 (\mathbb{H}^m \times \mathbb{H}^n,-\Delta - \lambda)$, but even the details of the partition of the full Martin boundary (\ref{part}) depend on $\lambda$.
\end{example}

\subsection{Boundaries and Singularities}\label{sec:applsing}

As an application of the theory developed above we get a full potential and Martin theory on uniform domains and minimal hypersurfaces.

\subsubsection{Uniform Spaces}\label{sec:appluniform}

Using the hyperbolic unfoldings presented in \autoref{sec:hypunf}, the results of the preceding settings can be easily translated to the setting of uniform spaces. In the statement of a boundary Harnack inequality, the sets $\N^\delta_i(\xi)$ from \thmref{lemma:den}, defined relative to the quasi-hyperbolic metric, behave similarly to concentric balls in the original metric.

\begin{proposition}[boundary Harnack inequality and Martin theory on uniform spaces]% 4.3
\label{cor:mtud} Let $\E=\E^0+V$ be a Schrödinger form satisfying a strong barrier condition on a bounded uniform space $Y$ with bounded geometry conditions as described in \autoref{sec:operatoronunifspaces}. Then there is a constant $C\ge 1$ depending only on universal constants such that for any two $\E$-harmonic functions $u$, $v >0$ on $\N^\delta_i(\xi)$ (defined in the quasi-hyperbolic metric) for some $\xi\in\del Y$, both $\L$-vanishing\footnote{The definition of $\L$-vanishing on $\N^\delta_i(\xi)$ is the same as in \thmref{Lvan} and \thmref{prop:LvanM} where we merely had an ideal boundary.} on $\N^\delta_i(\xi)$,
\[
\frac{u(x)}{v(x)} \le C \frac{u(y)}{v(y)} \text{ for any two points } x, y \in \N^\delta_{i+1}(\xi)\,.
\]
In particular, the topological and the Martin boundary are homeomorphic and every Martin boundary point is  minimal:
$\del(Y,d) \cong \delM^0(Y,\L) \cong \delM(Y,\L)$.
\end{proposition}
\begin{proof}
Using the hyperbolic unfolding from \autoref{sec:hypunf}, the boundary Harnack inequality follows from the \thmnameref{cor:hbhi}
applied to the $\L'$-harmonic functions $\dd^\frac{N-2}2u$ and $\dd^\frac{N-2}2v$ which are directly seen to be $\L'$-vanishing on $\N^\delta_i(\xi)$. Note that our more general notion of $\L$-vanishing carries over, while the functions $u$ and $v$ do not necessarily vanish towards the boundary in a classical sense. The additional factors $\dd^\frac{N-2}2$ cancel and we obtain the assertion.

Using \thmref{thm:grh}, \thmref{cor:gmr} and the transformation behavior of harmonic functions and Martin functions, we further conclude that
\[
    \del (Y,d) \cong  \delG (Y,d') \cong \delM (Y,\L')\cong  \delM (Y,\L)
\]
and every Martin boundary point is minimal.
\end{proof}

Boundary Harnack inequalities for harmonic functions on Lipschitz domains have first been investigated by Kemper in \cite{Kem72}.
While we followed the approach using hyperbolic unfoldings due to Ancona \cite{Anc87,Anc90} combined with geometric work by Bonk--Heinonen--Koskela \cite{BHK01}, there is also a more direct approach to boundary Harnack inequalities for the Laplacian on uniform domains initiated by Aikawa \cite{Aik01}, which does not require a strong barrier condition (corresponding to coercivity of $\L'$). More recently, Aikawa's ideas have been generalized to Dirichlet forms on metric spaces under various additional conditions \cite{LS14,BM19}.

\subsubsection{Minimal Hypersurfaces}\label{sec:applmini}

In the case of a singular minimal hypersurface $H$, we get similar results as in the case of uniform domains, with the important difference that the singular set $\Sigma$ becomes the Martin boundary. Since $H\setminus \Sigma$ is a manifold, one might choose to either apply the general theory or the simplifications in \autoref{rem:smooth}. \cite{Loh19} follows the second track.

\begin{theorem}[Boundary Harnack Inequality]\cite[Th.\,1 and 3]{Loh19} Let $\L$ be an \emph{\si-adapted} operator on $H\setminus \Sigma$, which implies our conditions in \autoref{sec:operatoronunifspaces} including a strong barrier. Then there exists a constant $C(H,\L) >1$ such that for any $\xi \in \widehat{\Sigma}$ and any two solutions
$u$, $v >0$ of $\L\, f= 0$ on $H \setminus \Sigma$ that are both $\L$-vanishing along $\mathcal{N}^\delta_i(\xi)\cap \widehat{\Sigma}$, we have the boundary Harnack inequality
\begin{equation*}%\label{fhepq1}
u(x)/v(x) \le C \cdot  u(y)/v(y) \mm{ \emph{for all} }x,\, y \in \mathcal{N}^\delta_{i+1}(\xi).
\end{equation*}
\begin{itemize}
  \item The identity map on $H \setminus \Sigma$ extends to a homeomorphism  between $\widehat{H}$ and the \emph{\textbf{Martin compactification}} $\overline{H \setminus \Sigma}_M$.
  \item All Martin boundary points are minimal:
  $\p^0_M (H \setminus \Sigma,\L) \cong \p_M(H \setminus \Sigma,\L)$.
\end{itemize}
In particular, $\widehat{\Sigma}$ and the minimal Martin boundary $\p^0_M (H \setminus \Sigma,\L)$ are homeomorphic.
\end{theorem}

Here, $\widehat{H}$ and $\widehat{\Sigma}$ denote the one-point compactifications of $H$ and $\Sigma$ in the non-compact case of Euclidean hypersurfaces with the following extra condition: we always add the point $\infty$
to $\Sigma$, even when $\Sigma$ is a compact subset of $H \subset \R^{n+1}$.

Minimal hypersurfaces have additional exploitable structure at the boundary in the form of (not necessarily unique) tangent cones in singular points. Any such tangent cone $C$ is again area-minimizing.
Their interplay with potential theory has the following consequences, see \cite[Th.\,3]{Loh19} and \cite[Th.\,3 and 4]{Loh21}:
\begin{itemize}
  \item If $\L_H$ is \si-adapted on $H$, $u>0$ solves $\L_H \, f = 0$ on $H \setminus \Sigma$ and $u$ is $\L_H$-vanishing along $B \cap \Sigma$ for some ball $B \subset H$ around a point $p \in \Sigma$, then the induced operator $\L_C$ on any of tangent cone $C$ of $H$ in $p$ is also \si-adapted. Moreover, any solution $v>0$ induced by $u$ under blow-up $\L_C$-vanishes along the total singular set $\Sigma_C$ of $C$.
  \item  The Martin theory of $\L_C$ on $C$ shows that there is exactly one solution $v>0$ of $\L_C\, f= 0$ with minimal growth towards $\Sigma_C$, up to multiples, namely the Martin function at $\infty\in \widehat{\Sigma_C}$. Moreover, $\L_C$  reproduces under scalings of $C$. Thus, we get a separation of variables: $v = \psi_C(\omega) \cdot r^{\alpha_C}$, $(\omega,r) \in  \p (B_1(0)  \cap C \setminus \Sigma_C) \times \R^{>0}$, for some function $\psi_C$ on $\p B_1(0)  \cap C \setminus \Sigma_C$ and $\alpha_C <0$.
  \item  This permits an asymptotic analysis of the solution $u$ on  $H \setminus \Sigma$ near $p$ from the family of Martin theories on $H \setminus \Sigma$ and on its tangent cones.
\end{itemize}

\paragraph{Acknowledgements:}
We thank the anonymous referee for many useful comments and suggestions that helped us to improve the exposition.
We acknowledge support from the Open Access Publication Fund of the University of Münster.

\paragraph{Conflict of Interest:} The authors state no conflict of interest.

\phantomsection\addcontentsline{toc}{section}{\refname}
\bibliographystyle{aomalpha2}
\bibliography{pgrcd}

\providecommand{\bysame}{\leavevmode\hbox to3em{\hrulefill}\thinspace}
\providecommand{\noopsort}[1]{}
\providecommand{\mr}[1]{\href{http://www.ams.org/mathscinet-getitem?mr=#1}{MR~#1}}
\providecommand{\zbl}[1]{\href{http://www.zentralblatt-math.org/zmath/en/search/?q=an:#1}{Zbl~#1}}
\providecommand{\jfm}[1]{\href{http://www.emis.de/cgi-bin/JFM-item?#1}{JFM~#1}}
\providecommand{\arxiv}[1]{\href{https://arxiv.org/abs/#1}{arXiv~#1}}
\providecommand{\doi}[1]{\href{https://doi.org/#1}{DOI~#1}}
\providecommand{\MR}{\relax\ifhmode\unskip\space\fi MR }
\providecommand{\MRhref}[2]{%
  \href{http://www.ams.org/mathscinet-getitem?mr=#1}{#2}
}
\providecommand{\href}[2]{#2}
\begin{thebibliography}{GMR15}

\bibitem[Aik01]{Aik01}
\bgroup\scshape{}H.~Aikawa\egroup{}, Boundary {{Harnack}} principle and
  {{Martin}} boundary for a uniform domain,  \emph{J. Math. Soc. Japan}
  \textbf{53} (2001), 119--145. \mr{1800526}.

\bibitem[Aik12]{Aik12}
\bgroup\scshape{}H.~Aikawa\egroup{}, Potential analysis on nonsmooth
  domains\textemdash{{Martin}} boundary and boundary {{Harnack}} principle,  in
  \emph{Complex Analysis and Potential Theory}, \emph{{{CRM Proc}}. {{Lecture
  Notes}}} \textbf{55}, {Amer. Math. Soc., Providence, RI}, 2012, pp.~235--253.
  \mr{2986906}.

\bibitem[AMR97]{AMR97}
\bgroup\scshape{}S.~Albeverio\egroup{}, \bgroup\scshape{}Z.-M. Ma\egroup{}, and
  \bgroup\scshape{}M.~R{\"o}ckner\egroup{}, Partitions of unity in {{Sobolev}}
  spaces over infinite-dimensional state spaces,  \emph{J. Funct. Anal.}
  \textbf{143} (1997), 247--268. \mr{1428125}.  \doi{10.1006/jfan.1996.2968}.

\bibitem[AH18]{AH18}
\bgroup\scshape{}L.~Ambrosio\egroup{} and \bgroup\scshape{}S.~Honda\egroup{},
  Local spectral convergence in {${\rm RCD}^*(K,N)$} spaces,  \emph{Nonlinear
  Anal.} \textbf{177} (2018), 1--23. \mr{3865185}.
  \doi{10.1016/j.na.2017.04.003}.

\bibitem[Anc90]{Anc90}
\bgroup\scshape{}A.~Ancona\egroup{}, Th\'eorie du potentiel sur les graphes et
  les vari\'et\'es,  in \emph{\'Ecole d'\'et\'e de {{Probabilit\'es}} de
  {{Saint-Flour XVIII}}\textemdash 1988}, \emph{Lecture {{Notes}} in {{Math}}.}
  \textbf{1427}, {Springer, Berlin}, 1990, pp.~1--112. \mr{1100282}.
  \doi{10.1007/BFb0103041}.

\bibitem[Anc86]{Anc86}
\bgroup\scshape{}A.~Ancona\egroup{}, On strong barriers and an inequality of
  {{Hardy}} for domains in {{$\mathbb{R}^{n}$}},  \emph{J. London Math. Soc.
  (2)} \textbf{34} (1986), 274--290. \mr{856511}.
  \doi{10.1112/jlms/s2-34.2.274}.

\bibitem[Anc87]{Anc87}
\bgroup\scshape{}A.~Ancona\egroup{}, Negatively curved manifolds, elliptic
  operators, and the {{Martin}} boundary,  \emph{Ann. of Math.} \textbf{125}
  (1987), 495--536. \mr{890161}.  \doi{10.2307/1971409}.

\bibitem[Anc88]{Anc88}
\bgroup\scshape{}A.~Ancona\egroup{}, Positive harmonic functions and
  hyperbolicity,  in \emph{Potential {{Theory Surveys}} and {{Problems}}}
  (\bgroup\scshape{}J.~Kr{\'a}l\egroup{}, \bgroup\scshape{}J.~Luke{\v
  s}\egroup{}, \bgroup\scshape{}I.~Netuka\egroup{}, and
  \bgroup\scshape{}J.~Vesel{\'y}\egroup{}, eds.), \textbf{1344}, {Springer
  Berlin Heidelberg}, {Berlin, Heidelberg}, 1988, pp.~1--23.
  \doi{10.1007/BFb0103341}.

\bibitem[Anc12]{Anc12}
\bgroup\scshape{}A.~Ancona\egroup{}, On positive harmonic functions in cones
  and cylinders,  \emph{Rev. Mat. Iberoam.} \textbf{28} (2012), 201--230.
  \mr{2904138}.  \doi{10.4171/RMI/674}.

\bibitem[AS85]{AS85}
\bgroup\scshape{}M.~T. Anderson\egroup{} and
  \bgroup\scshape{}R.~Schoen\egroup{}, Positive harmonic functions on complete
  manifolds of negative curvature,  \emph{Ann. of Math.} \textbf{121} (1985),
  429--461. \mr{794369}.  \doi{10.2307/1971181}.

\bibitem[BM19]{BM19}
\bgroup\scshape{}M.~T. Barlow\egroup{} and
  \bgroup\scshape{}M.~Murugan\egroup{}, Boundary {{Harnack}} principle and
  elliptic {{Harnack}} inequality,  \emph{J. Math. Soc. Japan} \textbf{71}
  (2019), 383--412. \mr{3943443}.  \doi{10.2969/jmsj/77057705}.

\bibitem[Bau66]{Bau66}
\bgroup\scshape{}H.~Bauer\egroup{}, \emph{{Harmonische R\"aume und ihre
  Potentialtheorie}}, \emph{{Lecture Notes in Mathematics, Vol. 22}},
  {Springer, Berlin-New York}, 1966. \mr{0210916}.

\bibitem[Bir94]{Bir94}
\bgroup\scshape{}M.~Biroli\egroup{}, The {{Wiener}} test for
  {{Poincar\'e-Dirichlet}} forms,  in \emph{Classical and Modern Potential
  Theory and Applications ({{Chateau}} de {{Bonas}}, 1993)}, \emph{{{NATO
  Adv}}. {{Sci}}. {{Inst}}. {{Ser}}. {{C}}: {{Math}}. {{Phys}}. {{Sci}}.}
  \textbf{430}, {Kluwer Acad. Publ., Dordrecht}, 1994, pp.~93--104.
  \mr{1321609}.

\bibitem[BM95]{BM95}
\bgroup\scshape{}M.~Biroli\egroup{} and \bgroup\scshape{}U.~Mosco\egroup{}, A
  {{Saint-Venant}} type principle for {{Dirichlet}} forms on discontinuous
  media,  \emph{Ann. Mat. Pura Appl. (4)} \textbf{169} (1995), 125--181.
  \mr{1378473}.  \doi{10.1007/BF01759352}.

\bibitem[BM07]{BM07}
\bgroup\scshape{}M.~Biroli\egroup{} and \bgroup\scshape{}S.~Marchi\egroup{},
  Wiener criterion at the boundary related to {\emph{p}}-homogeneous strongly
  local {{Dirichlet}} forms,  \emph{Matematiche (Catania)} \textbf{62} (2007),
  37--52. \mr{2401177}.

\bibitem[BB11]{BB11}
\bgroup\scshape{}A.~Bj{\"o}rn\egroup{} and
  \bgroup\scshape{}J.~Bj{\"o}rn\egroup{}, \emph{Nonlinear Potential Theory on
  Metric Spaces}, \emph{{{EMS Tracts}} in {{Mathematics}}} \textbf{17},
  {European Mathematical Society, Z\"urich}, 2011. \mr{2867756}.
  \doi{10.4171/099}.

\bibitem[BB18]{BB18}
\bgroup\scshape{}A.~Bj{\"o}rn\egroup{} and
  \bgroup\scshape{}J.~Bj{\"o}rn\egroup{}, Local and semilocal {{Poincar\'e}}
  inequalities on metric spaces,  \emph{J. Math. Pures Appl. (9)} \textbf{119}
  (2018), 158--192. \mr{3862146}.  \doi{10.1016/j.matpur.2018.05.005}.
  \arxiv{1703.00752}.

\bibitem[BMS01]{BMS01}
\bgroup\scshape{}J.~Bj{\"o}rn\egroup{}, \bgroup\scshape{}P.~MacManus\egroup{},
  and \bgroup\scshape{}N.~Shanmugalingam\egroup{}, Fat sets and pointwise
  boundary estimates for \$p\$-harmonic functions in metric spaces,  \emph{J.
  Anal. Math.} \textbf{85} (2001), 339--369. \mr{1869615}.
  \doi{10.1007/BF02788087}.

\bibitem[BB07]{BB07}
\bgroup\scshape{}S.~Blach{\`e}re\egroup{} and
  \bgroup\scshape{}S.~Brofferio\egroup{}, Internal diffusion limited
  aggregation on discrete groups having exponential growth,  \emph{Probab.
  Theory Related Fields} \textbf{137} (2007), 323--343. \mr{2278460}.
  \doi{10.1007/s00440-006-0009-2}.  \arxiv{math/0507582}.

\bibitem[BHM11]{BHM11}
\bgroup\scshape{}S.~Blach{\`e}re\egroup{},
  \bgroup\scshape{}P.~Ha{\"i}ssinsky\egroup{}, and
  \bgroup\scshape{}P.~Mathieu\egroup{}, Harmonic measures versus quasiconformal
  measures for hyperbolic groups,  \emph{Ann. Sci. \'Ec. Norm. Sup\'er. (4)}
  \textbf{44} (2011), 683--721. \mr{2919980}.  \doi{10.24033/asens.2153}.
  \arxiv{0806.3915}.

\bibitem[BH86]{BH86}
\bgroup\scshape{}J.~Bliedtner\egroup{} and \bgroup\scshape{}W.~Hansen\egroup{},
  \emph{Potential Theory}, \emph{Universitext}, {Springer, Berlin, Heidelberg},
  1986. \mr{850715}.  \doi{10.1007/978-3-642-71131-2}.

\bibitem[Bon96]{Bon96}
\bgroup\scshape{}M.~Bonk\egroup{}, Quasi-geodesic segments and {{Gromov}}
  hyperbolic spaces,  \emph{Geom. Dedicata} \textbf{62} (1996), 281--298.
  \mr{1406442}.  \doi{10.1007/BF00181569}.

\bibitem[BHK01]{BHK01}
\bgroup\scshape{}M.~Bonk\egroup{}, \bgroup\scshape{}J.~Heinonen\egroup{}, and
  \bgroup\scshape{}P.~Koskela\egroup{}, \emph{Uniformizing {{Gromov}}
  Hyperbolic Spaces}, \emph{Ast\'erisque} \textbf{270}, {Soci\'et\'e
  math\'ematique de France}, 2001. \mr{1829896}.

\bibitem[BJ06]{BJ06}
\bgroup\scshape{}A.~Borel\egroup{} and \bgroup\scshape{}L.~Ji\egroup{},
  \emph{Compactifications of Symmetric and Locally Symmetric Spaces},
  \emph{Mathematics: {{Theory}} \& {{Applications}}}, {Birkh\"auser, Boston},
  2006. \mr{2189882}.

\bibitem[Bre67]{Bre67}
\bgroup\scshape{}M.~Brelot\egroup{}, \emph{Lectures on Potential Theory},
  second ed., \emph{Tata {{Institute}} of {{Fundamental Research Lectures}} on
  {{Mathematics}}} \textbf{19}, {Tata Institute of Fundamental Research,
  Bombay}, 1967. \mr{0259146}.

\bibitem[BH99]{BH99}
\bgroup\scshape{}M.~R. Bridson\egroup{} and
  \bgroup\scshape{}A.~Haefliger\egroup{}, \emph{{Metric Spaces of Non-Positive
  Curvature}}, \emph{{Grundlehren der Mathematischen Wissenschaften}}
  \textbf{319}, {Springer, Berlin}, 1999. \mr{1744486}.
  \doi{10.1007/978-3-662-12494-9}.

\bibitem[BS07]{BS07}
\bgroup\scshape{}S.~Buyalo\egroup{} and \bgroup\scshape{}V.~Schroeder\egroup{},
  \emph{Elements of Asymptotic Geometry}, \emph{{{EMS Monographs}} in
  {{Mathematics}}}, {European Mathematical Society, Z\"urich}, 2007.
  \mr{2327160}.  \doi{10.4171/036}.

\bibitem[Cho69]{Cho69}
\bgroup\scshape{}G.~Choquet\egroup{}, \emph{Lectures on Analysis. {{Vol}}.
  {{I}}: {{Integration}} and Topological Vector Spaces}, {W. A. Benjamin, New
  York-Amsterdam}, 1969. \mr{0250011}.

\bibitem[CC72]{CC72}
\bgroup\scshape{}C.~Constantinescu\egroup{} and
  \bgroup\scshape{}A.~Cornea\egroup{}, \emph{Potential Theory on Harmonic
  Spaces}, {Springer-Verlag, New York-Heidelberg}, 1972. \mr{0419799}.

\bibitem[GO79]{GO79}
\bgroup\scshape{}F.~W. Gehring\egroup{} and \bgroup\scshape{}B.~G.
  Osgood\egroup{}, Uniform domains and the quasihyperbolic metric,  \emph{J.
  Analyse Math.} \textbf{36} (1979), 50--74. \mr{581801}.

\bibitem[Gd90]{Gd90}
\bgroup\scshape{}{\'E}.~Ghys\egroup{} and \bgroup\scshape{}P.~{de la
  Harpe}\egroup{} (eds.), \emph{{Sur les groupes hyperboliques d'apr\`es
  Mikhael Gromov}}, \emph{{Progress in Mathematics}} \textbf{83},
  {Birkh\"auser, Boston}, 1990. \mr{1086648}.

\bibitem[Gig15]{Gig15}
\bgroup\scshape{}N.~Gigli\egroup{}, On the differential structure of metric
  measure spaces and applications,  \emph{Mem. Amer. Math. Soc.} \textbf{236}
  (2015), vi+91. \mr{3381131}.  \doi{10.1090/memo/1113}.

\bibitem[GMR15]{GMR15}
\bgroup\scshape{}N.~Gigli\egroup{}, \bgroup\scshape{}A.~Mondino\egroup{}, and
  \bgroup\scshape{}T.~Rajala\egroup{}, Euclidean spaces as weak tangents of
  infinitesimally {{Hilbertian}} metric measure spaces with {{Ricci}} curvature
  bounded below,  \emph{J. Reine Angew. Math.} \textbf{705} (2015), 233--244.
  \mr{3377394}.  \doi{10.1515/crelle-2013-0052}.

\bibitem[GW93]{GW93}
\bgroup\scshape{}S.~Giulini\egroup{} and \bgroup\scshape{}W.~Woess\egroup{},
  The {{Martin}} compactification of the {{Cartesian}} product of two
  hyperbolic spaces,  \emph{J. Reine Angew. Math.} \textbf{444} (1993), 17--28.
  \mr{1241792}.  \doi{10.1515/crll.1993.444.17}.  Available at
  \url{http://gdz.sub.uni-goettingen.de/dms/resolveppn/?PPN=GDZPPN002211076}.

\bibitem[GH08]{GH08}
\bgroup\scshape{}A.~Grigor'yan\egroup{} and \bgroup\scshape{}J.~Hu\egroup{},
  Off-diagonal upper estimates for the heat kernel of the {{Dirichlet}} forms
  on metric spaces,  \emph{Invent. Math.} \textbf{174} (2008), 81--126.
  \mr{2430977}.  \doi{10.1007/s00222-008-0135-9}.

\bibitem[GH14]{GH14}
\bgroup\scshape{}A.~Grigor'yan\egroup{} and \bgroup\scshape{}J.~Hu\egroup{},
  Heat kernels and {{Green}} functions on metric measure spaces,  \emph{Canad.
  J. Math} \textbf{66} (2014), 641--699.

\bibitem[Gro87]{Gro87}
\bgroup\scshape{}M.~Gromov\egroup{}, Hyperbolic groups,  in \emph{Essays in
  Group Theory}, \emph{Math. {{Sci}}. {{Res}}. {{Inst}}. {{Publ}}.} \textbf{8},
  {Springer, New York}, 1987, pp.~75--263. \mr{919829}.

\bibitem[GS11]{GS11}
\bgroup\scshape{}P.~Gyrya\egroup{} and
  \bgroup\scshape{}L.~{Saloff-Coste}\egroup{}, Neumann and {{Dirichlet}} heat
  kernels in inner uniform domains,  \emph{Ast\'erisque} (2011). \mr{2807275}.

\bibitem[Hel69]{Hel69}
\bgroup\scshape{}L.~L. Helms\egroup{}, \emph{Introduction to Potential Theory},
  \emph{Pure and {{Applied Mathematics}}, {{Vol}}. {{XXII}}},
  {Wiley-Interscience A Division of John Wiley \& Sons, New
  York-London-Sydney}, 1969. \mr{0261018}.

\bibitem[Her62]{Her62}
\bgroup\scshape{}R.-M. Herv{\'e}\egroup{}, Recherches axiomatiques sur la
  th\'eorie des fonctions surharmoniques et du potentiel,  \emph{Ann. Inst.
  Fourier (Grenoble)} \textbf{12} (1962), 415--571. \mr{139756}.
  \doi{10.5802/aif.125}.

\bibitem[Kem72]{Kem72}
\bgroup\scshape{}J.~T. Kemper\egroup{}, A boundary {{Harnack}} principle for
  {{Lipschitz}} domains and the principle of positive singularities,
  \emph{Comm. Pure Appl. Math.} \textbf{25} (1972), 247--255. \mr{293114}.
  \doi{10.1002/cpa.3160250303}.

\bibitem[Leh17]{Leh17}
\bgroup\scshape{}J.~Lehrb{\"a}ck\egroup{}, Hardy inequalities and {{Assouad}}
  dimensions,  \emph{J. Anal. Math.} \textbf{131} (2017), 367--398.
  \mr{3631460}.  \doi{10.1007/s11854-017-0013-8}.

\bibitem[Lew88]{Lew88}
\bgroup\scshape{}J.~L. Lewis\egroup{}, Uniformly fat sets,  \emph{Trans. Amer.
  Math. Soc.} \textbf{308} (1988), 177--196. \mr{946438}.
  \doi{10.2307/2000957}.

\bibitem[LS14]{LS14}
\bgroup\scshape{}J.~Lierl\egroup{} and
  \bgroup\scshape{}L.~{Saloff-Coste}\egroup{}, Scale-invariant boundary
  {{Harnack}} principle in inner uniform domains,  \emph{Osaka J. Math.}
  \textbf{51} (2014), 619--656. \mr{3272609}.

\bibitem[LW65]{LW65}
\bgroup\scshape{}P.~A. Loeb\egroup{} and \bgroup\scshape{}B.~Walsh\egroup{},
  The equivalence of {{Harnack}}'s principle and {{Harnack}}'s inequality in
  the axiomatic system of {{Brelot}},  \emph{Ann. Inst. Fourier (Grenoble)}
  \textbf{15} (1965), 597--600. \mr{190360}.

\bibitem[Loh18]{Loh18}
\bgroup\scshape{}J.~Lohkamp\egroup{}, Hyperbolic unfoldings of minimal
  hypersurfaces,  \emph{Anal. Geom. Metr. Spaces} \textbf{6} (2018), 96--128.
  \mr{3849619}.  \doi{10.1515/agms-2018-0006}.

\bibitem[Loh20]{Loh19}
\bgroup\scshape{}J.~Lohkamp\egroup{}, Potential theory on minimal hypersurfaces
  {{I}}: {{Singularities}} as {{Martin}} boundaries,  \emph{Potential Analysis}
  \textbf{53} (2020), 1493--1528. \mr{4159389}.
  \doi{10.1007/s11118-019-09815-6}.

\bibitem[Loh21]{Loh21}
\bgroup\scshape{}J.~Lohkamp\egroup{}, Potential theory on minimal hypersurfaces
  {{II}}: {{Hardy}} structures and {{Schr\"odinger}} operators,
  \emph{Potential Analysis} \textbf{55} (2021), 563--602. \mr{4341062}.
  \doi{10.1007/s11118-020-09869-x}.

\bibitem[MR92]{MR92}
\bgroup\scshape{}Z.~M. Ma\egroup{} and
  \bgroup\scshape{}M.~R{\"o}ckner\egroup{}, \emph{Introduction to the Theory of
  (Nonsymmetric) {{Dirichlet}} Forms}, \emph{Universitext}, {Springer-Verlag,
  Berlin}, 1992. \mr{1214375}.  \doi{10.1007/978-3-642-77739-4}.

\bibitem[Mun00]{Mun00}
\bgroup\scshape{}J.~R. Munkres\egroup{}, \emph{Topology}, second ed., {Prentice
  Hall, Upper Saddle River, NJ}, 2000. \mr{3728284}.

\bibitem[Pin95]{Pin95}
\bgroup\scshape{}R.~G. Pinsky\egroup{}, \emph{Positive Harmonic Functions and
  Diffusion}, \emph{Cambridge {{Studies}} in {{Advanced Mathematics}}}
  \textbf{45}, {Cambridge University Press, Cambridge}, 1995. \mr{1326606}.
  \doi{10.1017/CBO9780511526244}.

\bibitem[Ste70]{Ste70}
\bgroup\scshape{}E.~M. Stein\egroup{}, \emph{Singular Integrals and
  Differentiability Properties of Functions}, \emph{Princeton {{Mathematical
  Series}}} \textbf{30}, {Princeton University Press, Princeton, N.J.}, 1970.
  \mr{0290095}.

\bibitem[Stu96]{Stu96}
\bgroup\scshape{}K.~T. Sturm\egroup{}, Analysis on local {{Dirichlet}} spaces.
  {{III}}. {{The}} parabolic {{Harnack}} inequality,  \emph{J. Math. Pures
  Appl. (9)} \textbf{75} (1996), 273--297. \mr{1387522}.

\end{thebibliography}

{\small
\noindent\textsc{Mathematisches Institut, Universität Münster, Münster, Germany}

\noindent\textit{E-mail addresses:} \email{m.kemper@uni-muenster.de}, \email{j.lohkamp@uni-muenster.de}
\end{document}